\documentclass[ejs]{imsart}
\RequirePackage[OT1]{fontenc}
\RequirePackage{amsthm,amsmath}
\RequirePackage[numbers]{natbib}
\RequirePackage[colorlinks,citecolor=blue,urlcolor=blue]{hyperref}
\usepackage{natbib}
\usepackage{psfrag, graphicx, euscript,amsfonts,latexsym,amssymb,latexsym}
\usepackage{epstopdf}

\renewcommand{\hat}{\widehat}
\renewcommand{\tilde}{\widetilde}
\newcommand{\rf}[1]{(\ref{eq:#1})}

\newcommand{\rE}{{\mathbb E}}

\newcommand{\rR}{\mathbb R}
\newcommand{\rd}{{\rm d}}
\def\mathbbone{{\mathchoice {\rm 1\mskip-4mu l} {\rm 1\mskip-4mu l}
{\rm 1\mskip-4.5mu l} {\rm 1\mskip-5mu l}}}
\newcommand{\ind}{{\mathbbone}} 
\newcommand{\pr} {\par \noindent{\bf Proof\,:~}}
\newcommand{\epr}{\hfill\hbox{\hskip 4pt
                \vrule width 5pt height 6pt depth 1.5pt}\vspace{0.5cm}\par}

\newcounter{assc}
\newcommand{\bcond}[1]{\refstepcounter{assc}\label{ass:#1} \begin{maliste}{\textbf{Condition K}}}
\newcommand{\econd}{\end{maliste}}

\newcommand{\bass}[1]{\refstepcounter{assc}\label{ass:#1} \begin{maliste}{\textbf{Assumption S}}}
\newcommand{\eass}{\end{maliste}}

\newtheorem{lemma}{Lemma}
\newtheorem{theorem}{Theorem}
\newtheorem{proposition}{Proposition}

\newtheorem{remark}{Remark}

\newtheorem{corollary}{Corollary}

\begin{document}

\begin{frontmatter}
\title{Optimal Rate of Direct Estimators in Systems of Ordinary Differential Equations Linear in Functions of the Parameters}
\runtitle{Estimation for ODE systems}


\author{\fnms{Itai} \snm{Dattner}\corref{}\ead[label=e1]{idattner@stat.haifa.ac.il}\thanksref{t1}}
\and
\author{\fnms{Chris A.J.} \snm{Klaassen}\ead[label=e2]{c.a.j.klaassen@uva.nl}\thanksref{t1}}

\thankstext{t1}{Research under STW-grant at EURANDOM}

\affiliation{University of Haifa and Korteweg-de Vries Institute for Mathematics, University of
Amsterdam}

\address{Department of Statistics,
University of Haifa,
199 Aba Khoushy Ave.
Mount Carmel, Haifa 3498838
, Israel.
\printead{e1}}

\address{Korteweg-de Vries Institute for Mathematics,
University of
Amsterdam,
P.O. Box 94248
1090 GE Amsterdam, The Netherlands.
\printead{e2}}

\runauthor{I. Dattner and C.A.J. Klaassen}

\begin{abstract}
Many processes in biology, chemistry, physics, medicine, and engineering are modeled by a system of differential equations. Such a system is usually characterized via unknown parameters and estimating their 'true' value is thus required. In this paper we focus on the quite common systems for which the derivatives of the states may be written as sums of products of a function of the states and a function of the parameters.

For such a system linear in functions of the unknown parameters we present a necessary and sufficient condition for identifiability of the parameters. We develop an estimation approach that bypasses the heavy computational burden of numerical integration and avoids the estimation of system states derivatives, drawbacks from which many classic estimation methods suffer. We also suggest an experimental design for which smoothing can be circumvented. The optimal rate of the proposed estimators, i.e., their $\sqrt n$-consistency, is proved and simulation results illustrate their excellent finite sample performance and compare it to other estimation approaches.
\end{abstract}


\tableofcontents
\end{frontmatter}

\section{Introduction}
Many processes in biology, chemistry, physics, medicine, and engineering are modeled by a system of differential equations. Parameter estimation for such systems is considered as the bottleneck in modeling dynamic processes and attracts some growing attention in recent statistical literature. In particular, new estimation methods are developed (e.g., \cite{ramsay2007parameter}, \cite{qi2010asymptotic}) or quite old techniques are rigorously analyzed (e.g., \cite{xue2010sieve}, \cite{gugushvili2012parametric}). Below we review other research as well. Most of it considers systems of ordinary differential equations (ODEs) of the form
\begin{equation}\label{eq:ode_model}
\bigg\{
\begin{array}{l}
x^{\prime}(t)=F(x(t);\nu),\ t\in[0,1],
\\
x(0)=\xi,
\end{array}
\end{equation}
where $x(t)$ takes values in $\rR^d,\, \xi$ in $\Xi\subset \rR^d,$
and $\nu \in N \subset {\mathbb R}^q$.
The seemingly more general nonautonomous system
\begin{equation*}
\bigg\{
\begin{array}{l}
\tilde{x}^{\prime}(t)=F(\tilde{x}(t),t;\nu),\ t\in[0,1],
\\
\tilde{x}(0)=\tilde{\xi},
\end{array}
\end{equation*}
may and will be reduced to the autonomous system (\ref{eq:ode_model}) by the simple substitution $x(t)=(\tilde{x}^T(t), t)^T,\ t\in[0,1],\xi=(\tilde{\xi}^T,0)^T.$

In many applications states and parameters can be separated in the sense that there exist measurable functions $g:\,\,{\mathbb R}^d \to {\mathbb R}^d \times {\mathbb R}^p$ and $h\,:\,N \to {\mathbb R}^p$ such that
\begin{equation}\label{eq:main_idea}
F(x(t);\nu)=g(x(t))h(\nu)
\end{equation}
holds. We write $\theta = h(\nu), \ \theta \in \Theta=h(N) \subset {\mathbb R}^p$, and call it the natural parameter, where $\nu$ is the parameter of interest.
\par
The class of ODEs (\ref{eq:main_idea}) is widely used in practice because of interpretability of the natural parameters as rate constants. In statistics a similar structure is popular; think of linear regression and e.g. Cox' proportional hazards model. The following list includes examples of systems in which the function $h$ is the identity, i.e., systems that are linear in the parameters: the Lotka-Volterra system in population dynamics (\cite{edelstein2005mathematical}); models describing HIV dynamics (\cite{nowak2000virus}, \cite{miao2008modeling}, \cite{miao2009differential}, \cite{wu2008parameter}, \cite{fang2011two}); models for the blood coagulation process (\cite{hockin2002model}); problems in chemistry (\cite{tjoa1991simultaneous}); gene regulatory networks (\cite{brewer2008fitting}); models describing the spread of infectious diseases (\cite{he2010plug},\cite{hooker2011parameterizing}); calcium measurements analysis (\cite{tank1995quantitative}); pharmacokinetic models (\cite{de2005comprehensive}). A well known example for the case where the system is not linear in the parameters but separability of the states and parameters is still possible, is the FitzHugh-Nagumo system in neurophysiology (\cite{fitzhugh1961impulses}, \cite{nagumo1962active}).
\par
The extensive list of applications above suggests that systems for which it is possible to separate the states from the parameters deserve special attention and treatment. However, current methods do not seem to exploit the full potential in such systems, both theoretically and practically. In the present study we attempt to do just this: in Section~\ref{sec:identifiability} we discuss identifiability in systems linear in the parameters; in Section~\ref{sec:methodology} we present a general estimation approach for the case where all trajectories of $x$ are observed. In Section~\ref{sec:asymptotic} we develop two estimators for different experimental setups and derive their $\sqrt n$-consistency, i.e., their optimal rate of convergence. In Section~\ref{sec:numerical} the application of the methods is demonstrated via simulations and a discussion is presented in Section~\ref{sec:discussion}. The proofs are relegated to the Appendix.
\section{Identifiability}\label{sec:identifiability}
A prerequisite for consistent estimation is that the parameter is
identifiable.  There are several concepts of identifiability (e.g., \cite{bellman1970structural}, \cite{cobelli1980parameter}, \cite{ljung1994global}, \cite{xia2003identifiability}; see also \cite{miao2011identifiability} and references therein).
We are concerned with {\it structural identifiability}, a property that depends on the mathematical structure of the model, but is not affected by the randomness of physical experiments. To be more specific, the identifiability criterion given in Proposition~\ref{pr:identifiability} below is
given in terms of a particular solution (i.e. set of trajectories) to the initial value problem. Clearly, a particular solution depends on elements of the experimental setup such as initial conditions and control parameters. Verifying the structural identifiability of a model is usually a difficult task that can be carried-out only in models of low dimensions (e.g., models describing HIV dynamics studied in \cite{miao2009differential}, \cite{wu2008parameter} and \cite{miao2011identifiability}).
\par
Exploiting linearity in the natural parameter $\theta$ we start with the following observation. By integration, (\ref{eq:ode_model}) and (\ref{eq:main_idea}) yield
the system of integral equations
\begin{equation}\label{eq:int_model}
x(t)=\xi + \int_0^t g(x(s))\, \rd s\,\theta,\ t\in[0,1].
\end{equation}
Given the values of $\xi$ and $\theta$ the solution of
(\ref{eq:ode_model}), (\ref{eq:main_idea}), and (\ref{eq:int_model}) is denoted by
\begin{equation*}\label{eq:solution}
x(t)=x(t; \theta, \xi),\ t\in[0,1].
\end{equation*}
In the present context identifiability means that knowledge of a solution $x(t),\,
t\in[0,1],$ for the system
(\ref{eq:ode_model}), (\ref{eq:main_idea}), and (\ref{eq:int_model}) yields the values of
the parameters $\xi$ and $\theta$. For $\xi=x(0)$ this is
obviously true, while identifiability for $\theta$ means that
\begin{equation}\label{eq:ident}
\theta^\prime\ne\theta\Rightarrow x(\cdot;\theta^\prime,\xi)\ne x(\cdot;\theta,\xi).
\end{equation}
From (\ref{eq:int_model}) we see that different values of $\theta$
may yield the same solution $x(t),\,t\in[0,1],$ if and only if
the $p$ columns of $g(x(t))$ are linearly dependent satisfying a nontrivial linear equation that is the same for Lebesgue almost all $t \in [0,1].$ 
This observation is generalized and formulated precisely in the proposition below. For its formulation we need some notation. Let $W$ be a symmetric $d \times d$-matrix of finite signed measures on $([0,1], \cal B)$ with $\cal B$ the sigma field of Borel sets, and let $x\,:\, [0,1] \to {\mathbb R}^d$ and $y\,:\, [0,1] \to {\mathbb R}^d$ be Borel measurable vector valued functions. We assume that $W$ is chosen in such a way that
\begin{equation}\label{innerproduct}
<x,y>_W = \int_0^1 x^T(t)\, \rd W(t)\, y(t)
\end{equation}
is a semidefinite inner product and
\begin{equation}\label{norm}
\parallel x \parallel_W = <x,x>^{1/2}_W
\end{equation}
is the corresponding seminorm. Note that in (\ref{innerproduct}) the integration with respect to $\rd W(t)$ includes $y(t).$ To clarify this notation we note the following. Let $\mu$ be a finite nonnegative measure on $([0,1], {\cal B})$ dominating all signed measures in the matrix $W$ (for example, the sum of the variations of the finite signed measures in $W$ will do). Denote by $w(\cdot)$ the $d \times d$-matrix of the Radon-Nikodym derivatives of the signed measures in $W$ with respect to $\mu.$ Now (\ref{innerproduct}) may be rewritten as
\begin{equation}
< x, y >_W = \int_0^1 x^T(t) w(t) y(t) \rd \mu(t).
\end{equation}
Note that the inner product from (\ref{norm}) introduces equivalence classes of measurable functions in that $x\,:\, [0,1] \to {\mathbb R}^d$ and $y\,:\, [0,1] \to {\mathbb R}^d$ are equivalent if and only if $\parallel x-y \parallel_W =0$ holds. We shall assume that $\parallel x \parallel_W =0$ implies that $x_i(t)=0$ for $W_{ii}$-almost all $t \in [0,1]$ and for $i=1, \dots, d.$ We shall assume also that 0 belongs to the support of $W_{ii}$ for $i=1, \dots, d.$ If $x\,:\, [0,1] \to {\mathbb R}^{d \times k}$ and $y\,:\, [0,1] \to {\mathbb R}^{d \times \ell}$ are measurable matrix valued functions, then $<x,y>_W$ will be interpreted as the $k \times \ell$ matrix of the inner products of the columns of $x$ and of $y.$ Denote the $d\times d$ identity matrix by $I_d$ and assume that the matrix
\begin{equation}\label{normingW}
\int_0^1 \rd W(t) = <I_d, I_d>_W = A_W
\end{equation}
is well-defined with finite entries and positive definite.
\begin{proposition}\label{pr:identifiability}
Let $\xi\in\Xi$ and $\theta\in\Theta$ with $\Theta$ an open subset of ${\mathbb R}^p.$ Let
$x(t)=x(t;\theta,\xi),\,t\in[0,1],$ satisfy the system
(\ref{eq:ode_model})--(\ref{eq:int_model}) and write
\begin{equation}\label{G}
G(t) = \int_0^t g(x(s))\,\rd s\,,\quad t \in
[0,1].
\end{equation}
Let W be a symmetric $d \times d$-matrix of signed measures as in (\ref{innerproduct}) satisfying (\ref{normingW}) and having the other properties mentioned above. Assume that the
$d \times p$- and $p \times p$-matrices
\begin{equation}\label{BC}
B_W = < I_d, G>_W, \quad C_W= <G, G>_W
\end{equation}
are well-defined with finite entries.

{\rm (i)} If $C_W$ is nonsingular then $A_W - B_W C_W^{-1}B_W^T$ is and
\begin{eqnarray}
\xi &=& \left(A_W - B_W C_W^{-1}B_W^T\right)^{-1} <I_d - G C_W^{-1}B_W^T ,x>_W, \label{xi} \\
\theta &=& C_W^{-1} \left( <G, x>_W -B_W^T \xi \right) \label{theta}
\end{eqnarray}
hold.

{\rm (ii)}
Conversely, if knowledge for all $i=1, \dots,d$ of $x_i(t)$ for $W_{ii}$-almost all $t\in[0,1]$ determines $\theta,$ then $C_W$
is nonsingular.
\end{proposition}
A proof of this proposition is given in Appendix A.1, but here we would like to note already that (\ref{xi}) and (\ref{theta}) follow from the fact that at its minimum 0 the derivatives of $\parallel x - \zeta - G \eta \parallel_W^2$ with respect to $\eta$ and $\zeta$ at $\theta$ and $\xi$ respectively, have to vanish. Note that $C_W$ is singular if and only if there exists a $p$-vector $\eta\ne 0$ with
\begin{equation}\label{singularB}
C_W \eta= <G, G\eta >_W =0,
\end{equation}
which implies
\begin{equation*}
\eta^T C_W \eta= \parallel G \eta\parallel_W^2=0.
\end{equation*}

Consequently, $G(t)\eta=0$ for Lebesgue almost all $t \in [0,1],$ if all $W_{ii}$ are equivalent to Lebesgue measure on the unit interval, and hence the $p$ columns of $g(x(t))$ in $\rR^d$ satisfy the same nontrivial linear relationship for Lebesgue almost all $t \in [0,1].$
Conversely, this linear relationship on $g(x(t)), t\in [0,1],$ implies
(\ref{singularB}) and hence the singularity of $C_W.$
\par
A careful examination of the proposition above reveals that uniqueness of the solution $x(t;\theta,\xi)$ (as a function of $t$) is not required for identifiability of the natural parameter. Note that uniqueness of solutions was previously assumed in \cite{qi2010asymptotic}, \cite{xue2010sieve}, and \cite{gugushvili2012sqrt} who dealt with the fully nonlinear case.
According to the Picard-Lindel\"of theorem existence and uniqueness of the solution $x(\cdot; \theta, \xi)$ in some neighborhood of 0 is guaranteed if the map $g(\cdot)$ is Lipschitz continuous; see also \cite[Chapter 2]{arnold1977}. However, consider for any
positive $\alpha \ne 1$ the (one dimensional) initial value
problem
\begin{equation*}
\bigg\{
\begin{array}{l}
x^{\prime}(t)=\alpha\theta x(t)^{(\alpha -1)/\alpha},\quad
t\in[0,1], \\
x(0)=\xi=0.
\end{array}
\end{equation*}
One may check that
\[x(t;\theta,0)=[\theta(t-\tau)]^\alpha \vee 0,\quad
t\in[0,1],\] is a solution for any $\tau\in[0,1)$. Hence, there
are infinitely many solutions for this initial value problem.
Nevertheless, the parameter $\theta$ is identifiable which may be
verified by calculating for $W$ the uniform distribution $C_W=\int_0^1 G^T(t) G(t)\,\rd
t=\theta^{2(\alpha -1)}(1-\tau)^{2\alpha +1}/(2\alpha +1)\ne 0$
for any $\tau\in[0,1)$.
\par
As for identifiability of the parameter of interest $\nu$ we note that part (i) of Proposition \ref{pr:identifiability} may be applied if the measurable parametrization function $h: N \to \Theta$ is injective, namely
\begin{equation}\label{nu}
\nu = h^{-1}\left( C_W^{-1} \left( <G, x>_W -B_W^T \xi \right) \right).
\end{equation}
If the natural parameter is not identifiable, the parameter of interest might be. However, we will not study this rather complicated situation here.
\section{Methodological approach}\label{sec:methodology}
In practice, the values of $\xi, \theta,$ and $\nu$ are unknown and
one usually observes $x(t; \theta,\xi),\ \theta=h(\nu),$ with noise and at certain time points
only. We denote the $n$ observations by
\begin{equation}\label{eq:observations}
Y(t_i)=x(t_i;\theta,\xi) +\varepsilon(t_i)\,,\quad i=1,\dots,n,\
0\leq t_1 \leq \dots \leq t_n \leq 1,
\end{equation}
where $\varepsilon(t_i)$ is the unobserved $d$-dimensional column
vector of measurement errors at time $t_i.$ This experimental setup is common (e.g., \cite{ramsay2007parameter} and \cite{gugushvili2012sqrt}), and many methods for estimating parameters in this context have been developed. For an extensive survey of recent developments in parameter estimation and
structure identification of biochemical and genomic systems, see
\cite{chou2009recent}. Since the list of estimation methods is exhaustive, a detailed review is not feasible, thus we will focus on the two most relevant techniques: the first is the nonlinear least squares (NLS) method that motivates our study, while the second is the two-step approach which we adopt.
\par
The classical nonlinear least squares method aims at minimizing over $\eta\in\Theta$ and $\zeta\in\Xi$ the function
\begin{equation*}
\frac{1}{n}\sum_{i=1}^n\parallel
Y(t_i)- x(t_i;\eta,\zeta)\parallel^2,
\end{equation*}
where $\parallel \cdot \parallel$ denotes the standard Euclidean norm.
Unless an exact solution $x(t_i;\eta,\zeta)$ is at hand, it is approximated via numerical integration, and the minimization of the criterion function is
carried-out by searching the parameter space for the global minimum.
Statistical properties of this method are studied in
\cite{xue2010sieve} for the situation that $\xi\in\Xi$ is
known. However, \cite{voit2004decoupling} demonstrate that
the need to repeat numerical integration multiple times might
increase the computational time for numerical integration up to
$95\%$ of the total computational time required for a gradient
based optimization method (even in low dimensional systems).

In order to bypass the burden of numerical integration, several {\it collocation} estimation methods were developed, such as the {\it two-step} technique (e.g., \cite{Bellman197126}, \cite{varah1982spline}) and {\it generalized profiling} (\cite{ramsay2007parameter}). The generalized profiling method is asymptotically efficient (\cite{qi2010asymptotic}) provided the distribution of the measurement errors is known, and can handle a variety of problems (\cite{hooker2011parameterizing}, \cite{xun2013parameter}). On the other hand, the two-step approach, although requiring the choice of some smoothing parameter, is relatively more straightforward to apply. Thus, a two-step method can serve as a preliminary step in the parameter estimation task, to be followed by applying more complex methods such as generalized profiling. This type of estimation strategy was successfully demonstrated in \cite{vujacicetal2014ode} for fully observed systems, and in \cite{dattner2015modelbased} for the partially observed case. 

The classical two-step approach works as follows. The observations are first smoothed, which results in an estimator $\hat{x}_n(\cdot)$ for the solution
$x(\cdot;\theta,\xi)$ of the system, and by differentiation in the
estimator $\hat{x}_n^{\prime}(\cdot)$ for
$x^\prime(\cdot;\theta,\xi).$ Then the estimator for $\theta$ is
the minimizer $\hat{\theta}_n$ over $\eta\in\Theta$ of the smooth
criterion function
\begin{equation}\label{eq:sme}
\int_0^1\parallel\hat{x}_n^{\prime}(t)-F(\hat{x}_n(t);\eta)\parallel^2
w(t)\, \rd t,
\end{equation}
where $w$ is an appropriate weight function. By estimating the "true" trajectories of the system and their derivatives, the two-step approach bypasses the need to integrate the system numerically and as a result, the parameter estimates can be computed extremely fast (\cite{brunel2008parameter}, \cite{liang2008parameter}). Under regularity conditions
\cite{gugushvili2012sqrt} show that this "smooth and match" estimator
(SME) $\hat{\theta}_n$ has the $\sqrt{n}$-rate of convergence to
$\theta.$ This is an example of the use of nonparametric "plug-in"
or substitution estimators (see \cite{goldstein1992optimal} and
\cite{bickel2003nonparametric}). When the system is linear in the parameters, (\ref{eq:sme}) can be minimized straightforwardly, as noted in \cite{brewer2008fitting}, \cite{fang2011two} and \cite{gugushvili2012sqrt}. However, their methods are based on estimates of derivatives, and it is
well known (see \cite{voit2000computational} and
\cite{chou2009recent}) that estimating derivatives from
noisy and sparse data may be rather inaccurate. Indeed, this problem attracted some attention (\cite{hall2013quick}, \cite{brunel2014parametric}). 
The methodology developed in the present paper is a two-step approach that does not require the estimation of derivatives. 
 Moreover, we also pay attention to estimation of the initial value $\xi.$ 

Let $\hat{x}_n(t), t\in[0,1],$ be an estimator of $x(t; \theta,\xi)$ based on the observations (\ref{eq:observations}). In view of (\ref{eq:int_model}) and in analogy to (\ref{eq:sme}) it makes sense to estimate the parameters $\theta$ and $\xi$ by minimizing
\begin{equation}\label{criterion}
\left|\left|
\hat{x}_n(t)-\zeta-\int_0^tg(\hat{x}_n(s))\,\rd s\,\eta
\right|\right|_{W_n}^2
\end{equation}
over $\eta$ and $\zeta,$ where $W_n$ is an appropriate $d \times d$-matrix of signed measures on $([0,1], \cal B)$ as in Proposition \ref{pr:identifiability}. Denote
\begin{eqnarray}\label{Ghat}
\lefteqn{\hat{G}_n(t) = \int_0^t g(\hat{x}_n(s))\,\rd s\,,\quad t \in
[0,1],\nonumber} \\
&& A_n = <I_d, I_d>_{W_n}, \quad \hat{B}_n = <I_d, \hat{G}_n>_{W_n}, \quad
\hat{C}_n = <\hat{G}_n, \hat{G}_n>_{W_n}.
\end{eqnarray}
Minimizing the criterion function (\ref{criterion}) with respect
to $\zeta$ and $\eta$ results in the direct estimators (cf. (\ref{xi}) and (\ref{theta}))
\begin{eqnarray}
\hat{\xi}_n &=& \left(A_n - \hat{B}_n \hat{C}_n^{-1}
\hat{B}_n^T\right)^{-1} <I_d - \hat{G}_n \hat{C}_n^{-1}
\hat{B}_n^T, \hat{x}_n >_{W_n}, \label{xihat} \\
\hat{\theta}_n &=& \hat{C}_n^{-1} \left(< \hat{G}_n, \hat{x}_n>_{W_n} - \hat{B}_n^T \hat{\xi}_n \right). \label{thetahat}
\end{eqnarray}
Note that these estimators are well-defined only if the inverse matrices in
(\ref{xihat}) and (\ref{thetahat}) exist. In case the initial
value $\xi$ is known, (\ref{thetahat}) may be used with
$\hat{\xi}_n$ replaced by $\xi.$

In order to estimate the parameter of interest $\nu$ we choose a distance function $d_n(\cdot, \cdot)$ on ${\mathbb R}^p$ and we choose $\hat{\nu}_n$ in such a way that
\begin{equation}\label{nuhat}
d_n\left( h\left(\hat{\nu}_n\right), \hat{\theta}_n \right) \leq \inf_{\nu \in N} d_n \left(h(\nu), \hat{\theta}_n \right) + \frac 1n
\end{equation}
holds. Of course, if the infimum is attained, we choose $\hat{\nu}_n$ as the minimizer.


The idea of an integral-based estimation approach as in (\ref{criterion}) appeared already in \cite{himmelblau1967determination}. These authors chose the $d \times d$-matrix $W_n$ to be a diagonal matrix with each diagonal element a weight function putting all its mass at the observation times $t_1, \dots, t_n.$ They proposed three specific weight functions, namely equal weights at all time points and two data dependent weight functions. These choices, with all their mass at the observation times, allow these authors to skip the smoothing step and to use $Y(t_i)$ instead of our $\hat{x}(t_i).$ This has the disadvantage that they had to consider multiple versions of their design, which they called runs, in order to obtain a good performance of their estimators. However, they did not derive statistical properties of their estimators. Their method is referred to in the chemical engineering literature as the 'direct integral method' and some papers revisited this idea (\cite{yermakova1982direct}, \cite{vajda1986direct}, and \cite{font1997testing}). In the next section, we introduce two modifications of the "direct integral method". These "modified integral methods" yield estimators with such desired statistical properties as consistency and the parametric $\sqrt n$ rate of convergence.
Still, the resulting estimators will not be statistically efficient. By a one step Newton-Raphson type of modification they can be turned into estimators equivalent to least squares estimators, and into efficient estimators when the distribution of the measurement errors is known, and even into semiparametrically efficient estimators when the distribution of the measurement errors is unknown. These modifications are under study (see e.g., \cite{dattner2015accelerated}). A possible way to apply our "modified integral methods" to general ODE systems, which are not necessarily linear in functions of the parameter, is under study as well.
\section{Asymptotic properties}\label{sec:asymptotic}
We start with some general asymptotic results for the estimation approach defined above and then we discuss two specific experimental set-ups.
Comparing our estimators (\ref{xihat}) and (\ref{thetahat}) to
(\ref{xi}) and (\ref{theta}) we see that they are consistent if
$g(\cdot)$ is continuous and $\hat{x}_n(\cdot)$ is a consistent
estimator of $x(\cdot)$ in an appropriate sense. Indeed, with the
notation $\parallel x \parallel_\infty = \sup_{t\in
[0,1]}\parallel x(t) \parallel$ we have the following result.
\begin{theorem}\label{th:const}
Let the model be defined by \rf{ode_model}--\rf{int_model} with
the map $g:\rR^d\rightarrow\rR^d\times\rR^p$ continuous. Fix
$\xi\in\Xi$ and $\theta\in\Theta$ and let
$x(\cdot)=x(\cdot;\theta, \xi)$ exist and be bounded on $[0,1],$
so
\begin{equation*}\label{boundednorm}
\parallel x \parallel_\infty  < \infty.
\end{equation*}
Let W be a symmetric $d \times d$-matrix of signed measures as in (\ref{innerproduct}) satisfying the conditions of Proposition \ref{pr:identifiability}. Furthermore, let the matrix $C_W$ from (\ref{BC}) be nonsingular, which implies that $\theta$ is identifiable via (\ref{theta}). Finally, let $\hat{x}_n(\cdot)$ be a consistent estimator of $x(\cdot)=x(\cdot;\theta, \xi)$ in the supnorm, i.e.,
\begin{equation}\label{supnorm}
\parallel \hat{x}_n - x \parallel_\infty \stackrel
P{\to} 0.
\end{equation}
If the sequence of matrices $W_n$ converges weakly to $W$ in the sense that the elements of $W_n$ converge weakly to the corresponding elements of $W,$ then the estimators $\hat{\xi}_n$ and $\hat{\theta}_n$ as presented in
(\ref{xihat}) and (\ref{thetahat}) are asymptotically well-defined and consistent, i.e.,
\[(\hat{\theta}_n, \hat{\xi}_n)\stackrel{P}{\to}(\theta,\xi)\]
holds as $n\to \infty.$ Moreover, if $\theta=h(\nu)$ holds, $d_n(x,y)/\parallel x-y \parallel$ are bounded away from 0 and infinity for all $x, y \in {\mathbb R}^d$ with $x \neq y,$ and $h^{-1}(\cdot)$ is continuous, then $\hat{\nu}_n$ as defined via (\ref{nuhat}) is asymptotically consistent as well.
\end{theorem}
Consequently, we have consistency of our estimators at all values
of the parameters for which the conditions are satisfied. In view of \rf{int_model}, $x(\cdot)$ is bounded if the map
$g(\cdot)$ is bounded.
\par
Note that if the system is not linear in its parameters then the
criterion function as in (\ref{eq:sme}) cannot be solved directly and
one needs to search the parameter space for the minimum. This
procedure requires that the criterion function separates the
parameter space well (cf. equation (3.9) in
\cite{gugushvili2012sqrt}). In our case this condition is
immediately satisfied.
\par
In order to get consistency at a certain rate we need stronger
conditions on $g(\cdot)$ and the estimator $\hat{x}_n(\cdot).$
\begin{theorem}\label{th:rootnconsistency}
Let the model be defined by \rf{ode_model}--\rf{int_model} with
the map $g:\rR^d\rightarrow\rR^d\times\rR^p$ twice continuously
differentiable. Fix $\xi\in\Xi$ and $\theta\in\Theta$ and let
$x(\cdot)=x(\cdot;\theta, \xi)$ exist and be bounded on $[0,1].$
Assume that $\theta$ is identifiable. Let the $d \times d$-matrices $W_n$ and $W$ be as defined in (\ref{innerproduct}) satisfying the conditions of Proposition 1. Let $\hat{x}_n(\cdot)$ be an
estimator of $x(\cdot)=x(\cdot;\theta, \xi)$ with
\begin{equation}\label{supnormc}
\parallel \hat{x}_n \parallel_\infty  = O_p(1),\quad
\parallel \rE\hat{x}_n -x\parallel_\infty \ = O(c_n),\quad c_n
\downarrow 0,
\end{equation}
and
\begin{equation}\label{dn}
\parallel \hat{x}_n - \rE\hat{x}_n \parallel^2_{W_n} = O_p(d_n), \quad d_n \downarrow 0.
\end{equation}
Assume that for every differentiable function $f: [0,1]\to \rR^d$ with bounded derivatives
\begin{equation}\label{orderweakconvergence}
<I_d, f>_{W_n} - <I_d, f>_W = O\left(w_n \right),\quad w_n \downarrow 0,
\end{equation}
holds. If for every bounded measurable function $f: [0,1]\to \rR,$
each component $\hat{x}_{n,j}(\cdot),\, j=1,\dots,d,$ of
$\hat{x}_n(\cdot),$ and all $h=1, \dots,d$
\begin{equation}\label{variance}
\int_0^1 {\rm var}\left(\int_0^t f(s)\, \hat{x}_{n,j}(s)\,\rd
s\right)\,\rd W_{n,hh}(t)=O\left(v_n^2 \right),\quad v_n \downarrow 0,
\end{equation}
and
\begin{equation}\label{variance2}
{\rm var}\left(\int_0^1 f(s)\, \hat{x}_{n,j}(s)\,\rd W_{n,hj}(s) \right)
=O\left(v_n^2 \right),\quad v_n \downarrow 0,
\end{equation}
hold, then estimators $\hat{\xi}_n$ and $\hat{\theta}_n$ as
defined in (\ref{xihat}) and (\ref{thetahat}) are consistent to
the following order
\begin{equation}\label{orderconsistency}
(\hat{\theta}_n -\theta, \hat{\xi}_n -\xi) = O_p\left( c_n + d_n +
v_n + w_n \right)
\end{equation}
as $n\to \infty.$ Furthermore, if
$g:\rR^d\rightarrow\rR^d\times\rR^p$ is twice differentiable and
all second derivatives of all components of $g(\cdot)$ are
bounded, then the condition $\parallel \hat{x}_n \parallel_\infty = O_p(1)$ is not needed in order to obtain
(\ref{orderconsistency}).
Moreover, if $\theta=h(\nu)$ holds, $d_n(x,y)/\parallel x-y \parallel$ are bounded away from 0 and infinity for all $x, y \in {\mathbb R}^d$ with $x \neq y,$ and $h^{-1}(\cdot)$ is Lipschitz continuous, then $\hat{\nu}_n$ as defined via (\ref{nuhat}) is asymptotically consistent to the order $O_p\left( c_n + d_n + v_n +w_n \right)$ as well.
\end{theorem}

Clearly with $c_n + d_n + v_n + w_n =O(n^{-1/2})$ this Theorem presents
sufficient conditions for the fastest possible rate, which means $\sqrt n$-consistency. In the next
subsection we present an estimator satisfying these conditions.
\subsection{Smooth estimator of solution ODE}\label{sec:smooth}

Our estimators ${\hat \theta}_n$ and ${\hat \xi}_n$ are defined by
(\ref{Ghat})--(\ref{thetahat}) and are based on an estimator
${\hat x}_n(\cdot)$ of the solution $x(\cdot)$ of the ODE system
(\ref{eq:ode_model})--(\ref{eq:int_model}). Clearly the quality of
the estimators ${\hat \theta}_n$ and ${\hat \xi}_n$ depends on the
properties of the estimator ${\hat x}_n(\cdot),$ as is illustrated
by the conditions of Theorem \ref{th:rootnconsistency}. Since the
classical kernel estimators are inconsistent at the boundaries of
the interval $[0,1]$, they do not satisfy these conditions.
Consequently, we need other estimators of $x(\cdot).$

Our choice here is to use a local polynomial type of estimator.
Under the assumption that all components of the solution
$x(\cdot)$ are $C^\alpha$-functions for some real $\alpha \geq 1,$ we
will approximate them by polynomials of degree
$\ell=\lfloor\alpha\rfloor$. This works as follows; cf.
\cite[Section 1.6]{tsybakov2009introduction}. For a given point $t_i\,,\,
i=1,\dots,n,$ and for $t$ sufficiently close to $t_i$ the
$d$-vector $x(t_i)$ equals approximately
\begin{eqnarray*}\label{LPestimator}
x(t_i)\approx x(t)+x^\prime(t)(t_i-t)+\cdots+x^{(\ell)}(t)
\frac{(t_i-t)^\ell}{\ell!}
=\nu(t)U\Big(\frac{t_i-t}{b}\Big),\nonumber
\\
U(u)=\Big(1,u,u^2/(2!),...,u^\ell/(\ell!)\Big)^T,\quad u\in \mathbb R,\\
\nu(t)=\left(x(t),x^\prime(t)b,x^{\prime\prime}(t)b^2,...,x^{(\ell)}(t)b^\ell\right),\quad
t \in \mathbb R, \nonumber
\end{eqnarray*}
where $b=b_n>0$ is a bandwidth, the $(\ell+1)$-vector $U(u)$ is a
column vector, and $\nu(t)$ is a $d \times (\ell+1)$-matrix. Let
$K(\cdot)$ be some appropriate kernel function and define
\[\hat{\nu}_n(t)=\arg\min_{\nu\in\rR^{d \times (\ell+1)}}
\sum_{i=1}^n \Big\{Y(t_i)-\nu U\Big(\frac{t_i-t}{b}\Big)\Big\}^T
\Big\{Y(t_i)-\nu U\Big(\frac{t_i-t}{b}\Big)\Big\}
K\Big(\frac{t_i-t}{b}\Big).\] The local polynomial estimator of
order $\ell$ of $x(t)$ is the first column of the $d\times (\ell
+1)$-matrix $\hat{\nu}_n(t)$, i.e.,
$\hat{x}_n(t)=\hat{\nu}_n(t)U(0).$ For a fixed $t$ this estimator
is just a weighted least squares estimator (\cite[Section
1.6]{tsybakov2009introduction}) and it may be written as the linear estimator
\begin{eqnarray}\label{eq:est_x}
\hat{x}_n(t)=\sum_{i=1}^n V_{n,i}(t)Y(t_i)
\end{eqnarray}
with
\begin{eqnarray*}\label{eq:W}
V_{n,i}(t)&=&\frac{1}{nb}U^T\Big(\frac{t_i-t}{b}\Big)B_n^{-1}(t)U(0)K\Big(\frac{t_i-t}{b}\Big),
\\\label{eq:B}
B_n(t)&=&\frac{1}{nb}\sum_{i=1}^nU\Big(\frac{t_i-t}{b}\Big)U^T\Big(\frac{t_i-t}{b}\Big)K\Big(\frac{t_i-t}{b}\Big).
\end{eqnarray*}
The following conditions on the kernel $K$ will assure that the
matrix $B_n(t)$ is positive definite and the estimator \rf{est_x}
is unique.\\ \medskip

\bcond{Condition K}
\begin{itemize}
\item[(i)] The kernel $K$ is symmetric around zero and has compact support, which lies within $[-1,1]$.
\item[(ii)] The kernel $K$ is Lipschitz on $\rR$, i.e., there exists a finite constant $L_K$
with  $|K(x)-K(y)|\leq L_K|x-y|, \forall\  x,y\in\rR$.
\item[(iii)] There exist constants $K_{\min}>0$, $\delta>0,$ and $K_{\max}<\infty$ with
$K_{\min}\ind_{[|u|\leq\delta]}\leq |K(u)|\leq K_{\max},
 \forall\
u\in\rR$.
\item[(iv)] The bandwidth $b=b_n$ satisfies $b_n\downarrow 0$ and $nb_n\rightarrow\infty$ as $n\rightarrow\infty$.
\end{itemize}
\econd
Conditions (i) and (iv) above are typical assumptions in
kernel estimation. The Lipschitz property in (ii) is needed when
deriving upper bounds for the risk of the estimator with respect
to the supremum norm. The lower bound for the kernel function in
(iii) is needed to assure that the matrix $B_n(t)$ is positive
definite.

Local polynomial estimators are consistent and "automatically"
correct for the boundaries. We note that some types of boundary
kernel estimators have bias and variance that are of the same
order. However, usually they have a complicated form and are not
easy to implement (see \cite{cheng1997automatic} for a
discussion on this problem). The following theorem assures us that
estimating $x(\cdot)$ by a local polynomial estimator fulfills the
requirements of Theorem~\ref{th:rootnconsistency}. A careful
choice of the bandwidth $b=b_n$ will result in a $\sqrt{n}$-rate
for the estimators $\hat{\theta}_n$ and $\hat\xi_n$.
\begin{theorem}\label{cor:local}
Let the model be defined by \rf{ode_model}--\rf{int_model} with
the map $g:\rR^d\rightarrow\rR^d\times\rR^p$ twice differentiable.
Fix $\xi\in\Xi$ and $\theta\in\Theta$ and let
$x(\cdot)=x(\cdot;\theta, \xi)$ exist. Suppose that for any
$j=1,...,d$ the component $x_j(t;\theta,\xi)$ is a
$C^\alpha$-function of $t$ on the interval $[0,1]$ for some real
$\alpha \geq 1$. Assume that $\theta$ is identifiable.

Let the observations be given by (\ref{eq:observations}) with
$t_i=i/n$, $i=1,...,n$. Assume that $\varepsilon_j(t_i),\
i=1,...,n,\ j=1,...,d,$ are i.i.d. with mean 0 and finite variance
$\sigma_\varepsilon^2.$ Let the estimator $\hat{x}_n(\cdot)$ for
$x(\cdot)$ be given in~\rf{est_x} with $\ell=\lfloor\alpha\rfloor$
and $b=b_n=n^{-1/(2\alpha)}.$

Let $W$ be a $d \times d$ matrix of signed measures as in Proposition \ref{pr:identifiability}. Let the estimators $\hat{\xi}_n$ and
$\hat{\theta}_n$ be defined in (\ref{xihat}) and (\ref{thetahat}) with the $d \times d$ matrix $W_n$ of signed measures satisfying
\begin{equation}\label{Wrootn1}
\sup_{1\leq h,j \leq d}\sup_{B \in {\cal B}([0,1])} \left| W_{n,hj}(B) - W_{hj}(B) \right| = O\left(\frac 1{\sqrt n} \right).
\end{equation}
Furthermore, let there exist a constant $C$ such that for any $h=1, \dots, d,$ for any interval $I_n$ of length $b_n\,,$ and for all $n$
\begin{equation}\label{Wrootn2}
W_{n,hh}\left(I_n \right) \leq C b_n
\end{equation}
holds. Under Assumption~$K$ the estimators $\hat{\theta}_n$ and
$\hat{\xi}_n$ are $\sqrt n$-consistent, i.e.,
\begin{equation}\label{smoothsrc}
\sqrt{n}(\hat{\theta}_n-\theta,\hat{\xi}_n-\xi) = O_p(1)
\end{equation}
holds, in the following cases:
\begin{enumerate}
\item
$\alpha \geq 3/2$ and $g(\cdot)$ has continuous second derivatives,
\item
$\alpha \geq 1$ and $g(\cdot)$ has bounded second derivatives.
\end{enumerate}
Moreover, if $\theta=h(\nu)$ holds, $d_n(x,y)/\parallel x-y \parallel$ are bounded away from 0 and infinity for all $x, y \in {\mathbb R}^d$ with $x \neq y,$ and $h^{-1}(\cdot)$ is Lipschitz continuous, then $\hat{\nu}_n$ as defined via (\ref{nuhat}) is $\sqrt n$-consistent as well.
\end{theorem}
%
%
%
%
%
Condition (\ref{Wrootn1}) states that the total variation distance between $W_n$ and $W$ should converge to 0 sufficiently fast. Note that (\ref{Wrootn2}) is satisfied if the $W_{n,hh}$ have bounded densities with respect to Lebesgue measure on $[0,1]$ or with respect to $1/n$ times counting measure on $i/n,\, i=1, \dots,n.$ Furthermore, note that for any $\theta\in\Theta$ and $\xi\in\Xi$ the
component $x_j(t;\theta,\xi)$ of the solution is a
$C^\alpha$-function in $t$ in a neighborhood of $0$ provided the
map $g$ is $C^\alpha$ in its argument (\cite[p. 52, Section
7.6, Corollary 4]{arnold1977}).
\par
Notice that \cite{gugushvili2012sqrt} study systems that are not necessarily
linear in the parameters. To prove $\sqrt n$-consistency of their
estimator they need Gaussianity or boundedness of the measurement
errors. Here just mean 0 and finite variance suffice.
\par

The method developed above is based on the preliminary step of
smoothing the observations. As a result, the performance of this
method is heavily based on the choice of the smoothing parameter.
This choice is not trivial in practice (see e.g.,
\cite{ramsay2007parameter}, \cite{qi2010asymptotic}
and \cite{gugushvili2012sqrt}), especially if one deals with a large
system and if the underlying system has "fast" and "slow"
components. In that case, using different bandwidths for different
components makes more sense. However, the proof of Theorem 3 will show that for $\alpha \geq 3/2$ the choice $b\approx n^{-1/3}$ always suffices.

\subsection{Step function estimator of solution ODE}\label{sec:step}
As mentioned above, choosing the smoothing parameter in practice may not be trivial. This problem can be avoided in situations like the following repeated measures model,
\begin{equation}\label{eq:observation_model_repeated}
Y^{(j)}(t_i)=x(t_i;\theta,\xi)+\varepsilon^{(j)}(t_i),\quad
j=1,\dots, J_i,\quad i=1,\dots,I,
\end{equation}
with $t_i=i/I,\ i=1,\dots,I.$ Hence, we observe $J_i$ repeated
measures of $x(t_i)$ for each time point $t_i,$ which means that
we have $n=\sum_{i=1}^I J_i$ observations in total. This is common
practice in many fields and therefore makes a quite reasonable
experimental setup.

Within this observation scheme it is natural to estimate $x(t_i)$
by
\begin{equation*}
\hat{x}_n(t_i)=\frac{1}{J_i}\sum_{j=1}^{J_i} Y^{(j)}(t_i)
\end{equation*}
and even to estimate $x(t)$ by
\begin{equation}\label{eq:xhatrepeated}
\hat{x}_n(t)=\frac{1}{J_i}\sum_{j=1}^{J_i} Y^{(j)}(t_i),\quad
(i-1)/I < t \leq i/I,\quad i=1,\dots,I,
\end{equation}
where we complete the definition of $\hat{x}_n(t)$ on $[0,1]$ by
$\hat{x}_n(0)= \hat{x}_n(t_1).$ This definition does not mean that
we intend to estimate the initial value $x(0;\theta, \xi)=\xi$ by
$\hat{x}_n(0).$ The estimator $\hat{x}_n(\cdot)$ is a preliminary
estimator of $x(\cdot;\theta, \xi)$ that will be used to construct
a more accurate estimator of $\xi$ than $\hat{x}_n(0).$ We choose
$W_n$ as in Theorem \ref{cor:local}. Again, our estimators ${\hat \theta}_n$ and ${\hat \xi}_n$ are defined by
(\ref{Ghat})--(\ref{thetahat}).

This estimator $\hat{\theta}_n$ with $\hat{\xi}_n$ replaced by $\xi$ equals the estimator based on the direct integral method of
\cite{himmelblau1967determination} with $J_i=J$ the number of runs and with $W_{ij}=1$ in (5)--(8) of \cite{himmelblau1967determination}, provided
the starting values for all runs are the same and known.
Both estimators $\hat{\theta}_n$ and $\hat{\xi}_n$ are $\sqrt n$-consistent if the number of time points
$I$ is of order $\sqrt n$ and for most time points $t_i$ the
sample size $J_i$ is of order $\sqrt n$ too. We formulate this
accurately in the following theorem.
\begin{theorem}\label{th:step}
Let the model be defined by \rf{ode_model}--\rf{int_model} with
the map $g:\rR^d\rightarrow\rR^d\times\rR^p$ twice differentiable.
Fix $\xi\in\Xi$ and $\theta\in\Theta$ and let
$x(\cdot)=x(\cdot;\theta, \xi)$ exist and be bounded on $[0,1].$
Assume that $\theta$ is identifiable. Let the observations be
given by \rf{observation_model_repeated} with $t_i=i/I,\
i=1,...,I$. Assume that $\varepsilon_h^{(j)}(t_i),\ h=1,\dots, d,\
j=1,\dots, J_i,\ i=1,\dots,I,$ are i.i.d. random variables with
zero expectation and finite variance $\sigma_\varepsilon^2$. Let $W_n$
satisfy (\ref{Wrootn1}) and (\ref{Wrootn2}) with $b_n$ replaced by $1/I.$
Let $\hat{x}_n(\cdot)$ be given by (\ref{eq:xhatrepeated}) and let
$\hat{\theta}_n$ and $\hat{\xi}_n$ be defined in (\ref{thetahat}) and (\ref{xihat}).
Furthermore, let the sample sizes satisfy
\begin{equation}\label{samplesizes}
\liminf_{n \to \infty} \frac I{\sqrt n} > 0,\quad \limsup_{n\to
\infty} \sum_{i=1}^I \frac 1{J_i} < \infty,\quad \sum_{i=1}^I J_i
=n.
\end{equation}
If the second derivatives of each component of $g(\cdot)$ are
continuous or bounded, then
\begin{equation}
\sqrt{n}(\hat{\theta}_n -\theta, \hat{\xi}_n -\xi)= O_p(1),\quad n
\to \infty,
\end{equation}
holds. Moreover, if $\theta=h(\nu)$ holds, $d_n(x,y)/\parallel x-y \parallel$ are bounded away from 0 and infinity for all $x, y \in {\mathbb R}^d$ with $x \neq y,$ and $h^{-1}(\cdot)$ is Lipschitz continuous, then $\hat{\nu}_n$ as defined via (\ref{nuhat}) is $\sqrt n$-consistent as well.
\end{theorem}
Note that ${\rm var}\left(\hat{x}_n(0)\right) =
\sigma_\varepsilon^2/J_1$ holds, and that estimating $\xi$
via $\hat{x}_n(0)$ would not yield the best possible rate, unless
$J_1$ is of exact order $n.$ Indeed, the $\sqrt n$-rate is
achievable by $\hat{\theta}_n$ using the
information from all $I$ time points.\\

\section{Simulation study}\label{sec:numerical}
In our simulation study we report on the finite sample properties of the smooth estimator of Section~\ref{sec:smooth}, and the step function estimator of Section~\ref{sec:step}. The smooth estimator is tested by comparing its performance to that of the derivative based two-step approach and of the generalized profiling estimator. This comparison is done for the same situations as have been used in the simulation studies for these estimators in literature. 
The study of the step function estimator is focused on understanding
the effect on the estimation accuracy of the number of repeated measures, as well as of different error distributions. In all simulations below, whenever the integral approach is applied, the initial values are considered as unknown and therefore are estimated as well.

\subsection{Smooth estimator}
Several researchers studied the problem of parameter estimation for the FitzHugh-Nagumo model (\cite{fitzhugh1961impulses},
\cite{nagumo1962active}) and therefore it is a good example to consider. This is a system with two states proposed as a simplification of the model presented in \cite{hodgkin1952quantitative} for
studying and simulating the animal nerve axon. Specifically, this model is used
in neurophysiology as an approximation of the observed spike
potential and takes the form
\begin{equation}\label{eq:FitzHugh}
\bigg\{
\begin{array}{l}
x_1^{\prime}(t)=\gamma(x_1(t)-x^3_1(t)+x_2(t)),
\\
x_2^{\prime}(t)=-(1/\gamma)(x_1(t)-\alpha+\beta x_2(t)).
\end{array}
\end{equation}
The voltage $x_1(t)$ moving across the cell membrane depends
on the recovery variable $x_2(t)$.
\par
This system was studied in \cite{liang2008parameter} who applied the derivative based method and in \cite{ramsay2007parameter} who used generalized profiling. We will compare the integral based approach to the results in the aforementioned papers. Note that the FitzHugh-Nagumo model was studied also by \cite{campbell2012smooth} who pointed out some difficulties in estimating the parameters for this ODE system.
\subsubsection{Comparison  with the derivative based method}
By setting $\nu=(\alpha,\beta,\gamma)^T$, the system (\ref{eq:FitzHugh}) takes the form (\ref{eq:main_idea}) with $\theta=h(\nu)=(\gamma,1/\gamma,\alpha/\gamma,\beta/\gamma)^T$ and the corresponding matrix $g$ is
\begin{eqnarray*}
g(x(t))=\left(
\begin{array}{cccc}
x_1(t)-x^3_1(t)+x_2(t)&0&0&0
\\
0&-x_1(t)&1&-x_2(t)
\end{array}
\right).
\end{eqnarray*}
While estimating parameters using a derivative based method does not require knowledge of the initial condition vector $\xi$, this is not the case with the integral based approach. Therefore we consider the initial values to be unknown and estimate them as well.
\par
The experimental setup follows that of \cite{liang2008parameter}. The true parameter vector is set to $(\alpha,\beta,\gamma)^T=(0.34,0.2,3)^T$ and the initial conditions to $\xi=(\xi_1,\xi_2)^T=(0,0.1)^T$. The two signals are first generated by solving the system at $0.1$ time units on the interval $[0,20]$ ($n=201$; note that the theory as developed for the time interval $[0,1]$ in the preceding sections is valid for any bounded interval $[0,T]$ as may be seen by scaling.) and then we add Gaussian measurement errors with zero mean and variances $\sigma^2_1, \sigma^2_2$ respectively. In particular, here we used local polynomial estimators of order $\ell=1$ for estimating the two components of $x(\cdot)$. The kernel function used for generating the local polynomial estimators was the same one as considered in \cite{liang2008parameter}, namely, $K(t)=3/4(1-t^2)\textbf{1}\{|t|\leq 1\}$, where $\textbf{1}\{\cdot\}$ stands for the indicator function. The last choice that has to be made before proceeding, is that of the bandwidth $b$. As pointed out in Remark~3 of \cite{liang2008parameter}, the bandwidth selection is critical in local polynomial regression. They used a bandwidth that under-smooths with respect to the optimal bandwidth for estimating $x(\cdot).$ Here we simply choose $b=n^{-1/3}$ (see the proof of Theorem~\ref{cor:local}).
\par
Once $\hat\theta_n$ is obtained, we can estimate $\nu$ using (\ref{nuhat}). To be more specific, we take $d_n$ to be the Mahalanobis distance:
\begin{eqnarray}\label{eq:funest}
d_n(x,y)=\sqrt{(x-y)^T\hat\Sigma_n^{-1}(x-y)},
\end{eqnarray}
where $\hat\Sigma_n$ is the estimated covariance matrix of $\hat\theta_n$. Given the observations model (\ref{eq:observations}), it is natural to define a bootstrap procedure for estimating $\hat\Sigma_n$ as follows (cf. \cite{hardle1988bootstrapping}).
Repeat $B$ times the following steps:
\begin{enumerate}
\item[{\rm (i)}] For each point $t_i$ generate residuals 
$\tilde\varepsilon_j(t_i) = Y_j(t_i)-\hat x_j(t_i)$. 
\item[{\rm (ii)}] Center the residuals: 
$\bar\varepsilon_j(t_i) = \tilde\varepsilon_j(t_i)- \frac{1}{n}\sum_{i=1}^n\tilde\varepsilon_j(t_i)$.
\item[{\rm (iii)}] Sample $n$ residuals (with replacement) from $\bar\varepsilon_j(t_1),...,\bar\varepsilon_j(t_n)$ to obtain the bootstrap residuals $\varepsilon^{\star}_j(t_1),...,\varepsilon^{\star}_j(t_n)$.
\item[{\rm (iv)}] Set $Y^{\star}_j(t_i) = \hat x_j(t_i) +\varepsilon^{\star}_j(t_i),\   i=1,...,n$.
\end{enumerate}
We then use the bootstrap sample $Y^{\star}_j(t_1),\dots,Y^{\star}_j(t_n)$ and apply the estimation procedure. Denote the estimator for the vector $\theta$ in the $b$th bootstrap sample by $\hat\theta_{n,b}^*$ and its corresponding average over the $B$ bootstrap samples by $\overline{\theta_n^*}$. Then we define
\[\hat\Sigma_n=\frac 1B\sum_{b=1}^B\Big\{(\hat\theta_{n,b}^*-\overline{\theta_n^*})(\hat\theta_{n,b}^*-\overline{\theta_n^*})^T\Big\}.
\]
Then we minimize $d_n \left(h(\nu), \hat{\theta}_n \right)$ over $\nu$ using a standard nonlinear optimization procedure (in this case, function \text{fminsearch} in Matlab). As an initial guess for the optimization step we take an arbitrary estimate for $\nu$ denoted by $\hat\nu_0=(\hat\alpha_0,\hat\beta_0,\hat\gamma_0)^T$. In this case we obtained it as follows. Let $\hat{\theta}_n^1,\hat{\theta}_n^2,\hat{\theta}_n^3,\hat{\theta}_n^4$ stand for the components of the vector of estimates $\hat{\theta}_n$. Then $\hat\nu_0=(\hat{\theta}_n^1\hat{\theta}_n^3,\hat{\theta}_n^1\hat{\theta}_n^4,\hat{\theta}_n^1)^T$.
\par
We conducted $M=500$ Monte Carlo simulations as in \cite{liang2008parameter}. We set $B=100$ for the bootstrap samples. The resulting empirical means and standard deviations of the integral approach are displayed in Table~\ref{tab:compare_liang_wu1} where $36$ different variance combinations are considered. The estimation results are substantially better uniformly over the experimental study, than those reported in Table~1 of \cite{liang2008parameter} for the derivative based approach. Furthermore, another measure of accuracy presented in the aforementioned paper is the average relative estimation error (ARE). The ARE of a real-valued parameter $a$ over the $M$ Monte Carlo simulations is defined as
\[ARE(a)=\frac{1}{M}\sum_{m=1}^M\frac{|\hat a_m-a|}{|a|}\times 100 \%,\]
where $\hat a_m$ is an estimator of $a$ in simulation $m$, and in our case $M=500$.
Table~\ref{tab:compare_liang_wu2} here presents the ARE of the integral based two-step approach, and corresponds to Table~2 of \cite{liang2008parameter}. For convenience, the results of Table~2 of \cite{liang2008parameter} are presented in Table~\ref{tab:compare_liang_wu2} as well under the title "Derivative" since their method is a derivative based two-step approach.
We see that the ARE's of the integral approach are substantially better, uniformly over the experimental study, than those of the derivative based approach.
\begin{table}
\caption{\label{tab:compare_liang_wu1}\textbf{Empirical means (standard deviations) for estimating the parameters of the FitzHugh-Nagumo system}. The parameters are estimated by (\ref{xihat})--(\ref{thetahat}) and (\ref{eq:funest}), with (\ref{eq:est_x}) a local polynomial estimator of order $\ell=1$ and bandwidth $b=n^{-1/3}$. Based on $500$ Monte Carlo simulations. The true parameter vector is  $(\alpha,\beta,\gamma)^T=(0.34,0.2,3)^T$. The two signals are first generated by solving the system at $0.1$ time units on the interval $[0,20]$ ($n=201$) and then adding Gaussian measurement errors with zero mean and variances $\sigma^2_1, \sigma^2_2$ respectively.}																			
\centering
\fbox{%
\begin{tabular}{llrrr}
\multicolumn{1}{l}{}&	
\multicolumn{1}{l}{}&																					
\multicolumn{3}{c}{Parameters }\\																					
\hline	
\multicolumn{1}{l}{$\sigma^2_1$}&				
\multicolumn{1}{l}{$\sigma^2_2$}&	
\multicolumn{1}{c}{$\alpha$}&																					
\multicolumn{1}{c}{$\beta$}&																					
\multicolumn{1}{c}{$\gamma$}\\		
\hline																
  0.050 &    0.050 &     0.339 (   0.004)&    0.200 (   0.022) &    3.005 (   0.033) \\ 
  &    0.060  &    0.340 (   0.005)&    0.203 (   0.024) &    3.004 (   0.039)  \\ 
  &    0.070  &    0.340 (   0.004)&    0.202 (   0.029) &    3.004 (   0.043)  \\ 
  &    0.080  &    0.340 (   0.005)&    0.204 (   0.030) &    3.004 (   0.046)  \\ 
  &    0.090  &    0.340 (   0.005)&    0.202 (   0.034) &    3.010 (   0.055)  \\ 
  &    0.100  &    0.340 (   0.005)&    0.206 (   0.038) &    3.006 (   0.059)  \\ 
 0.060 &    0.050 &     0.339 (   0.005)&    0.201 (   0.023) &    2.999 (   0.034) \\ 
  &    0.060  &    0.340 (   0.005)&    0.200 (   0.026) &    3.002 (   0.039)  \\ 
  &    0.070 &     0.340 (   0.005)&    0.202 (   0.030) &    3.005 (   0.047) \\ 
  &    0.080  &    0.340 (   0.005)&    0.201 (   0.033) &    3.008 (   0.049)  \\ 
 &    0.090 &     0.340 (   0.006)&    0.204 (   0.036) &    3.006 (   0.058) \\ 
  &    0.100  &    0.340 (   0.006)&    0.202 (   0.039) &    3.001 (   0.064)  \\ 
0.070 &    0.050 &     0.339 (   0.005)&    0.201 (   0.024) &    2.998 (   0.036) \\ 
  &    0.060  &    0.339 (   0.005)&    0.199 (   0.027) &    2.999 (   0.041)  \\ 
  &    0.070  &    0.339 (   0.006)&    0.200 (   0.029) &    3.001 (   0.047)  \\ 
  &    0.080  &    0.339 (   0.006)&    0.200 (   0.031) &    3.003 (   0.052)  \\ 
  &    0.090  &    0.340 (   0.006)&    0.204 (   0.037) &    3.003 (   0.055)  \\ 
  &    0.100  &    0.340 (   0.007)&    0.206 (   0.041) &    3.002 (   0.065)  \\ 
    0.080 &    0.050 &     0.338 (   0.006)&    0.201 (   0.024) &    3.001 (   0.038) \\ 
  &    0.060  &    0.339 (   0.006)&    0.201 (   0.028) &    2.999 (   0.042)  \\ 
  &    0.070  &    0.339 (   0.006)&    0.204 (   0.029) &    3.003 (   0.048)  \\ 
  &    0.080  &    0.339 (   0.006)&    0.199 (   0.035) &    2.997 (   0.052)  \\ 
  &    0.090  &    0.339 (   0.007)&    0.197 (   0.035) &    3.000 (   0.056)  \\ 
  &    0.100  &    0.340 (   0.007)&    0.206 (   0.041) &    2.994 (   0.064)  \\ 
    0.090 &    0.050 &     0.339 (   0.006)&    0.201 (   0.025) &    2.998 (   0.039) \\ 
  &    0.060  &    0.339 (   0.006)&    0.201 (   0.031) &    3.000 (   0.046)  \\ 
  &    0.070  &    0.339 (   0.007)&    0.201 (   0.031) &    3.000 (   0.051)  \\ 
  &    0.080  &    0.339 (   0.007)&    0.200 (   0.034) &    3.004 (   0.055)  \\ 
  &    0.090  &    0.339 (   0.007)&    0.201 (   0.038) &    2.997 (   0.061)  \\ 
  &    0.100  &    0.339 (   0.008)&    0.202 (   0.043) &    2.999 (   0.066)  \\ 
0.100 &    0.050 &     0.338 (   0.007)&    0.198 (   0.028) &    3.000 (   0.041) \\ 
  &    0.060  &    0.339 (   0.007)&    0.201 (   0.030) &    2.999 (   0.048)  \\ 
  &    0.070  &    0.339 (   0.007)&    0.202 (   0.034) &    2.995 (   0.051)  \\ 
  &    0.080  &    0.339 (   0.008)&    0.202 (   0.035) &    3.000 (   0.056)  \\ 
  &    0.090  &    0.338 (   0.008)&    0.202 (   0.041) &    2.998 (   0.060)  \\ 
  &    0.100  &    0.339 (   0.008)&    0.198 (   0.040) &    3.000 (   0.063)  
\end{tabular}    }
\end{table}						
\begin{table}
\caption{\label{tab:compare_liang_wu2}\textbf{Comparison of the empirical relative estimation error (ARE) in FitzHugh-Nagumo system, of the integral approach (Integral) with that of the derivative based two-step method (Derivative) of \cite{liang2008parameter}}. For the Integral estimator the parameters are estimated by (\ref{xihat})--(\ref{thetahat}) and (\ref{eq:funest}), with (\ref{eq:est_x}) a local polynomial estimator of order $\ell=1$ and bandwidth $b=n^{-1/3}$. Based on $500$ Monte Carlo simulations. The true parameter vector is  $(\alpha,\beta,\gamma)^T=(0.34,0.2,3)^T$. The two signals are first generated by solving the system at $0.1$ time units on the interval $[0,20]$ ($n=201$) and then adding Gaussian measurement errors with zero mean and variances $\sigma^2_1, \sigma^2_2$ respectively.}																			
\centering
\fbox{%
\begin{tabular}{llrrrrrr}
\multicolumn{1}{l}{}&	
\multicolumn{1}{l}{}&																					
\multicolumn{3}{c}{Integral }&
\multicolumn{3}{c}{Derivative }\\																						
\hline	
\multicolumn{1}{l}{$\sigma^2_1$}&				
\multicolumn{1}{l}{$\sigma^2_2$}&	
\multicolumn{1}{c}{$\alpha$}&																					
\multicolumn{1}{c}{$\beta$}&																					
\multicolumn{1}{c}{$\gamma$}&	
\multicolumn{1}{c}{$\alpha$}&																					
\multicolumn{1}{c}{$\beta$}&																					
\multicolumn{1}{c}{$\gamma$}\\		
\hline																
0.05	&	0.05	&	0.99	&	8.85	&	0.91	&	6.21	&	17.77	&	16.33	\\
	&	0.06	&	1.06	&	9.90	&	1.04	&	7.27	&	17.36	&	15.83	\\
	&	0.07	&	1.04	&	11.83	&	1.13	&	7.21	&	20.63	&	15.66	\\
	&	0.08	&	1.15	&	12.17	&	1.24	&	7.17	&	26.96	&	14.53	\\
	&	0.09	&	1.22	&	13.60	&	1.52	&	7.27	&	30.60	&	14.16	\\
	&	0.10	&	1.28	&	15.03	&	1.57	&	7.72	&	24.42	&	14.08	\\
0.06	&	0.05	&	1.13	&	9.26	&	0.90	&	6.70	&	16.66	&	18.38	\\
	&	0.06	&	1.13	&	10.29	&	1.04	&	7.33	&	18.00	&	17.76	\\
	&	0.07	&	1.24	&	11.94	&	1.23	&	6.06	&	20.85	&	17.27	\\
	&	0.08	&	1.26	&	13.14	&	1.34	&	5.75	&	26.67	&	16.97	\\
	&	0.09	&	1.38	&	14.45	&	1.53	&	7.32	&	22.79	&	16.55	\\
	&	0.10	&	1.51	&	15.33	&	1.68	&	7.90	&	29.71	&	16.07	\\
0.07	&	0.05	&	1.24	&	9.55	&	0.96	&	6.44	&	14.62	&	19.22	\\
	&	0.06	&	1.30	&	11.15	&	1.09	&	7.70	&	18.72	&	18.65	\\
	&	0.07	&	1.38	&	11.75	&	1.25	&	7.95	&	17.30	&	18.59	\\
	&	0.08	&	1.41	&	12.84	&	1.42	&	6.66	&	19.37	&	18.08	\\
	&	0.09	&	1.50	&	14.86	&	1.49	&	8.18	&	27.57	&	17.63	\\
	&	0.10	&	1.59	&	16.87	&	1.73	&	8.09	&	29.94	&	18.14	\\
0.08	&	0.05	&	1.36	&	9.80	&	1.00	&	6.28	&	16.41	&	20.94	\\
	&	0.06	&	1.51	&	10.95	&	1.10	&	6.90	&	21.51	&	20.14	\\
	&	0.07	&	1.53	&	11.90	&	1.30	&	7.33	&	18.55	&	20.07	\\
	&	0.08	&	1.47	&	13.85	&	1.40	&	7.95	&	21.39	&	20.23	\\
	&	0.09	&	1.65	&	14.31	&	1.52	&	7.78	&	25.05	&	18.61	\\
	&	0.10	&	1.70	&	16.51	&	1.73	&	7.75	&	30.93	&	18.86	\\
0.09	&	0.05	&	1.48	&	10.01	&	1.03	&	7.31	&	17.76	&	21.77	\\
	&	0.06	&	1.49	&	12.25	&	1.21	&	7.22	&	21.76	&	21.48	\\
	&	0.07	&	1.57	&	12.57	&	1.35	&	7.38	&	15.44	&	21.18	\\
	&	0.08	&	1.64	&	13.53	&	1.47	&	7.38	&	22.85	&	20.30	\\
	&	0.09	&	1.69	&	14.99	&	1.61	&	7.04	&	28.70	&	20.33	\\
	&	0.10	&	1.86	&	17.21	&	1.76	&	8.45	&	29.78	&	20.39	\\
0.10	&	0.05	&	1.61	&	11.17	&	1.11	&	6.42	&	18.89	&	22.68	\\
	&	0.06	&	1.62	&	12.04	&	1.29	&	6.78	&	19.33	&	21.87	\\
	&	0.07	&	1.77	&	13.54	&	1.37	&	6.62	&	22.09	&	21.79	\\
	&	0.08	&	1.79	&	13.85	&	1.47	&	7.80	&	23.20	&	22.12	\\
	&	0.09	&	1.84	&	15.87	&	1.60	&	8.30	&	24.40	&	20.85	\\
	&	0.10	&	1.92	&	15.91	&	1.68	&	8.57	&	26.50	&	20.99	\\
\end{tabular}    }
\end{table}						
\subsubsection{Comparison with generalized profiling}
The experimental setup here follows that of \cite{ramsay2007parameter}. In particular, they consider the following FitzHugh-Nagumo model
\begin{equation}\label{eq:FitzHugh2}
\bigg\{
\begin{array}{l}
x_1^{\prime}(t)=c(x_1(t)-x^3_1(t)/3+x_2(t)),
\\
x_2^{\prime}(t)=-(1/c)(x_1(t)-a+b x_2(t)).
\end{array}
\end{equation}
The true parameter vector is set to $\theta=(a,b,c)^T=(0.2,0.2,3)^T$ and the initial conditions to $\xi=(\xi_1,\xi_2)^T=(-1,1)^T$. The two signals are first generated by solving the system at $0.05$ time units on the interval $[0,20]$ ($n=401$) and then we add Gaussian measurement errors with zero mean and variances $\sigma^2_1=\sigma^2_2=0.5$.

The integral approach is executed as described above, the initial conditions are estimated as well. The estimation results, based on $500$ Monte Carlo simulations, are presented in Table~\ref{tab:compare_ramsay}. Also, in the table we present the results of the generalized profiling estimator that is chosen to be adapted to the Gaussianity of the measurement errors, as reported in Table~1 of \cite{ramsay2007parameter}. However, since it is not clear to us which initial guess was used there for the optimization over the parameter space, we also generated one experiment of our own. In particular, we first generate an initial guess in the parameter space that follows a Gaussian random vector with means the true parameters and a standard deviation of $0.5$ (variable \text{jitter} in the original code downloaded from the authors website). Then we start the $500$ Monte Carlo simulations using the same initial guess all over. The results are similar to those reported in \cite{ramsay2007parameter} except for the parameter $c$ for which the variability is higher here. We did not repeat the same experiment for other initial guesses since depending on the distance of the random guess from the true parameter vector, it could take the program about $90$ seconds to execute only one simulation out of the $500$ simulations required (using Intel(R) Core(TM) i7-4550U CPU @ 1.50GHz 2.10GHz 64-bit). In comparison, using the same hardware, one simulation of computing the integral estimator (including generating the $100$ bootstrap samples for estimating the covariance) takes about $17$ seconds to conclude. We note that when the system is linear in the parameters then there is no need for the bootstrap and the execution time of the integral estimator drops to less than $0.2$ seconds. Also, in calculating the total execution time for the integral estimator we exclude the time needed for constructing the matrix $W$, the weights of the local polynomials, since this matrix can be constructed {\it before} any observations are generated.

In summary, the estimated variance of the generalized profiling estimator is smaller than that of the two-step based integral approach. This is not surprising, since the generalized profiling estimator is asymptotically efficient as it has been chosen to be adapted to the Gaussianity of the measurement errors. However, the generalized profiling approach
involves an iterative optimization method (Gauss-Newton), which in turn, requires a good initial guess in the parameter space. Otherwise, the resulting estimates and execution time may be very bad. Thus, the integral approach may be used as a preliminary step in the estimation procedure, since it provides theoretical and practical guarantees that the resulting estimates are in the vicinity of the true parameter vector. Such a strategy may substantially improve the execution time of the generalized profiling approach even for systems of small dimensions (see for example Table~3 in \cite{vujacicetal2014ode}).
\begin{table}
\caption{\label{tab:compare_ramsay}\textbf{Comparison of the empirical means (standard deviation) in FitzHugh-Nagumo system, of the integral approach with that of the generalized profiling of \cite{ramsay2007parameter}}. For the Integral estimator the parameters are estimated by (\ref{xihat})--(\ref{thetahat}) and (\ref{eq:funest}), with (\ref{eq:est_x}) a local polynomial estimator of order $\ell=1$ and bandwidth $b=n^{-1/3}$. Based on $500$ Monte Carlo simulations. The true parameter vector is  $(a,b,c)^T=(0.2,0.2,3)^T$. The two signals are first generated by solving the system at $0.05$ time units on the interval $[0,20]$ ($n=401$) and then adding Gaussian measurement errors with zero mean and variances $\sigma^2_1=\sigma^2_2=0.5$.}																																				 
\centering
\fbox{%
\begin{tabular}{llll}
\multicolumn{1}{l}{}&																					
\multicolumn{1}{c}{$a$}&																					
\multicolumn{1}{c}{$b$}&																					
\multicolumn{1}{c}{$c$}\\		
\hline	
Generalized profiling (Ramsay et al. (2007))&0.2005 (0.0149)& 0.1984 (0.0643)& 2.9949 (0.0264) \\
Generalized profiling (here)&0.2003 (0.0166)&    0.1986 (0.0679)&    3.0010 (0.0795)\\
Integral estimator &   0.1906 (0.0307) &   0.1859 (0.0905) &   2.9249 (0.1216) \\ 
\end{tabular}    }
\end{table}						

\begin{figure}[h]
\[
\resizebox{300pt}{200pt}{
 \includegraphics{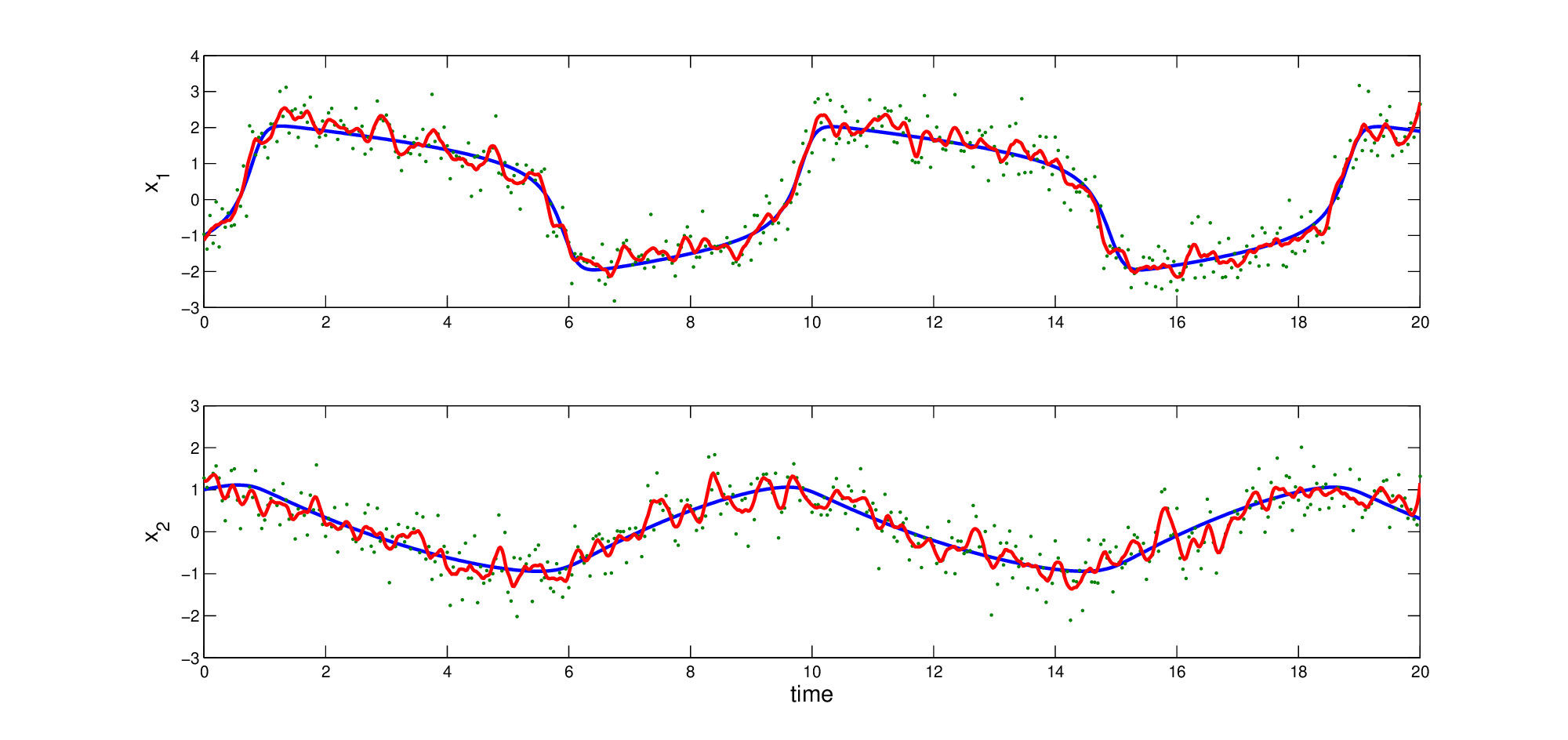}
}
\]
\caption{The FitzHugh-Nagumo system (setup as in Table \ref{tab:compare_ramsay}). \textbf{Top panel}. The blue
line corresponds to the true solution $x_1$; the dots correspond to one
realization of observations; the red line corresponds to the nonparametric estimator of $x_1$. \textbf{Bottom panel}. The same as the top but for $x_2$.
\label{fig:fitz}}
\end{figure}
\subsection{Step function estimator}
The goal of the following simulation study is merely to have a better understanding of the finite sample behavior of the step function estimator for different  repeated measures and noise scenarios. We consider the Lotka-Volterra system, a population dynamics model that describes evolution over time of the populations of two species, predators and their preys. In mathematical terms the
Lotka-Volterra model is described by a system consisting of two
equations and depending on the parameter
$\theta=(\theta_1,\theta_2,\theta_3,\theta_4)^T$. The system takes
the form
\begin{equation}\label{eq:lotka}
\bigg\{
\begin{array}{l}
x_1^{\prime}(t)=\theta_1x_1(t)-\theta_2x_1(t)x_2(t),
\\
x_2^{\prime}(t)=-\theta_3x_2(t)+\theta_4x_1(t)x_2(t).
\end{array}
\end{equation}
Here $x_1$ represents the size of the prey population and $x_2$ of
the predator population.

In the experiment we set the errors to be i.i.d. Gaussian or Laplace with zero mean and $\sigma^2_\varepsilon=0.5$ for both system states. In Tables~\ref{tab:normal}-\ref{tab:laplace} we present the empirical mean and standard deviation (in parenthesis) of the estimators for two different sets of parameters and initial values of the Lotka-Volterra system. Results are based on $5000$ Monte Carlo simulations. Both the rate constants $\theta$ and the initial values $\xi$ are estimated. In each simulation the data consist of $I=30$ noisy observations of $x_1$ and $x_2$
according to measurement error model \rf{observation_model_repeated}.
The samples were taken at $0.5$ time units on the interval $[0,T]$, $T=14.9$ for the first parameters setup and at $1$ time units on the interval $[0,T]$, $T=29.9$ for the second. At each time point, $J$ repeated measures were generated. 
Last two lines in each block correspond to the empirical mean and standard deviation (in parentheses) of the distribution of
$\{\frac{1}{T}\int_0^T \parallel x(t; \hat{\theta}_n,\hat{\xi}_n) - x(t;
\theta, \xi) \parallel^2 \rd t\}^{1/2}$ and $\parallel x(\cdot; \hat{\theta}_n,\hat{\xi}_n) - x(\cdot;
\theta, \xi) \parallel_\infty$ respectively.
The simulation results suggest that the finite sample behavior of the estimator is similar under both error distributions. Also, as expected, the estimation  accuracy grows with the number of repeated measures and is reasonable already when their number is relatively small.

\begin{table}
\caption{\label{tab:normal}\textbf{Gaussian error - empirical means (standard deviation) in Lotka-Volterra system based on $5000$ Monte Carlo simulations}. In each simulation the data consist of $I=30$ noisy observations of $x_1$ and $x_2$
according to measurement error model \rf{observation_model_repeated} with $\sigma^2_\varepsilon=0.5$. The samples were taken at $0.5$ time units on the interval $[0,14.9]$ for the first parameters setup and at $1$ time units on the interval $[0,29.9]$ for the second. At each time point, $J$ repeated measures were generated. The parameters are estimated by (\ref{xihat})--(\ref{thetahat}) with (\ref{eq:xhatrepeated}) a step function estimator. Last two lines in each block correspond to the empirical means (standard deviation)  of the distribution of
$\{\frac{1}{T}\int_0^T \parallel x(t; \hat{\theta}_n,\hat{\xi}_n) - x(t;
\theta, \xi) \parallel^2 \rd t\}^{1/2}$ and $\parallel x(\cdot; \hat{\theta}_n,\hat{\xi}_n) - x(\cdot;
\theta, \xi) \parallel_\infty$ respectively.}																					
\centering
\fbox{%
\begin{tabular}{llrrrr}
\hline	
\multicolumn{1}{l}{}&				
\multicolumn{1}{l}{Value}&																					
\multicolumn{1}{c}{$J=6$}&																					
\multicolumn{1}{c}{$J=10$}&																					
\multicolumn{1}{c}{$J=15$}&
\multicolumn{1}{c}{$J=30$}\\		
\hline																
$\xi_1$ &    1.000 &    1.089 (   0.143) &    1.085 (   0.116) &    1.085 (   0.093) &    1.083 (   0.065)\\
 $\xi_2$ &    0.500 &    0.446 (   0.096) &    0.441 (   0.075) &    0.438 (   0.061) &    0.436 (   0.043)\\
 $\theta_1$ &    0.500 &    0.468 (   0.077) &    0.473 (   0.061) &    0.474 (   0.050) &    0.477 (   0.035) \\
 $\theta_2$ &    0.500 &    0.473 (   0.075) &    0.477 (   0.060) &    0.479 (   0.048) &    0.480 (   0.034) \\
 $\theta_3$ &    0.500 &    0.500 (   0.073) &    0.501 (   0.057) &    0.501 (   0.047) &    0.501 (   0.033) \\
 $\theta_4$ &    0.500 &    0.508 (   0.073) &    0.508 (   0.057) &    0.509 (   0.047) &    0.509 (   0.033) \\
  &&    0.214 (   0.082) &    0.183 (   0.066) &    0.167 (   0.056) &    0.148 (   0.041) \\
  &&    0.344 (   0.142) &    0.290 (   0.111) &    0.264 (   0.094) &    0.233 (   0.069) \\ \hline
$\xi_1$ &    0.500 &    0.289 (   0.160) &    0.296 (   0.130) &    0.301 (   0.109) &    0.300 (   0.076)\\
 $\xi_2$ &    1.000 &    1.000 (   0.237) &    1.038 (   0.185) &    1.052 (   0.153) &    1.070 (   0.107)\\
 $\theta_1$ &    0.200 &    0.174 (   0.040) &    0.178 (   0.031) &    0.180 (   0.026) &    0.182 (   0.019) \\
 $\theta_2$ &    0.700 &    0.496 (   0.159) &    0.525 (   0.136) &    0.536 (   0.115) &    0.546 (   0.085) \\
 $\theta_3$ &    0.300 &    0.305 (   0.085) &    0.316 (   0.069) &    0.318 (   0.058) &    0.320 (   0.041) \\
 $\theta_4$ &    0.500 &    0.477 (   0.126) &    0.483 (   0.096) &    0.481 (   0.078) &    0.479 (   0.054) \\
  &&    0.443 (   0.249) &    0.385 (   0.189) &    0.341 (   0.166) &    0.284 (   0.137) \\
  &&    1.003 (   0.609) &    0.890 (   0.401) &    0.798 (   0.365) &    0.674 (   0.312) \\ \hline
\end{tabular}    }
\end{table}

\begin{table}
\caption{\label{tab:laplace}\textbf{Laplace error - empirical means (standard deviation) in Lotka-Volterra system based on $5000$ Monte Carlo simulations}. In each simulation the data consist of $I=30$ noisy observations of $x_1$ and $x_2$
according to measurement error model \rf{observation_model_repeated} with $\sigma^2_\varepsilon=0.5$. The samples were taken at $0.5$ time units on the interval $[0,14.9]$ for the first parameters setup and at $1$ time units on the interval $[0,29.9]$ for the second. At each time point, $J$ repeated measures were generated. The parameters are estimated by (\ref{xihat})--(\ref{thetahat}) with (\ref{eq:xhatrepeated}) a step function estimator. Last two lines in each block correspond to the empirical means (standard deviation)  of the distribution of
$\{\frac{1}{T}\int_0^T \parallel x(t; \hat{\theta}_n,\hat{\xi}_n) - x(t;
\theta, \xi) \parallel^2 \rd t\}^{1/2}$ and $\parallel x(\cdot; \hat{\theta}_n,\hat{\xi}_n) - x(\cdot;
\theta, \xi) \parallel_\infty$ respectively.}							
\centering
\fbox{%
\begin{tabular}{llrrrr}
\hline	
\multicolumn{1}{l}{}&				
\multicolumn{1}{l}{Value}&																					
\multicolumn{1}{c}{$J=6$}&																					
\multicolumn{1}{c}{$J=10$}&																					
\multicolumn{1}{c}{$J=15$}&
\multicolumn{1}{c}{$J=30$}\\																		
\hline
 $\xi_1$ &    1.000 &    1.087 (   0.147) &    1.088 (   0.112) &    1.086 (   0.092) &    1.085 (   0.065)\\
 $\xi_2$ &    0.500 &    0.444 (   0.097) &    0.440 (   0.074) &    0.439 (   0.062) &    0.436 (   0.044)\\
 $\theta_1$ &    0.500 &    0.470 (   0.079) &    0.472 (   0.059) &    0.474 (   0.049) &    0.476 (   0.035) \\
 $\theta_2$ &    0.500 &    0.475 (   0.076) &    0.477 (   0.057) &    0.478 (   0.048) &    0.480 (   0.034) \\
 $\theta_3$ &    0.500 &    0.499 (   0.073) &    0.501 (   0.056) &    0.501 (   0.046) &    0.502 (   0.032) \\
 $\theta_4$ &    0.500 &    0.507 (   0.073) &    0.509 (   0.057) &    0.508 (   0.046) &    0.510 (   0.033) \\
  &&    0.217 (   0.085) &    0.184 (   0.064) &    0.166 (   0.056) &    0.148 (   0.041) \\
  &&    0.348 (   0.149) &    0.292 (   0.109) &    0.262 (   0.092) &    0.233 (   0.069) \\ \hline
$\xi_1$ &    0.500 &    0.287 (   0.161) &    0.297 (   0.132) &    0.297 (   0.109) &    0.300 (   0.077)\\
 $\xi_2$ &    1.000 &    0.996 (   0.240) &    1.037 (   0.184) &    1.053 (   0.150) &    1.069 (   0.106)\\
 $\theta_1$ &    0.200 &    0.173 (   0.039) &    0.178 (   0.031) &    0.181 (   0.026) &    0.181 (   0.019) \\
 $\theta_2$ &    0.700 &    0.491 (   0.155) &    0.523 (   0.134) &    0.538 (   0.115) &    0.545 (   0.086) \\
 $\theta_3$ &    0.300 &    0.305 (   0.085) &    0.314 (   0.067) &    0.316 (   0.056) &    0.320 (   0.040) \\
 $\theta_4$ &    0.500 &    0.480 (   0.125) &    0.481 (   0.096) &    0.478 (   0.075) &    0.479 (   0.054) \\
  &&    0.446 (   0.380) &    0.379 (   0.184) &    0.339 (   0.167) &    0.281 (   0.135) \\
  &&    1.007 (   1.221) &    0.879 (   0.399) &    0.796 (   0.365) &    0.667 (   0.305) \\ \hline
\end{tabular}  }
\end{table}		

\section{Discussion}\label{sec:discussion}
Systems of ordinary differential equations are widely used by scientists for modeling real life phenomena. In this paper we studied systems for which separability of the states and parameters is possible, or more specifically, systems that are linear in functions of the parameters. Such systems are spread over diverse fields such as population dynamics, neurophysiology, HIV dynamics, blood coagulation, chemistry, gene regulatory networks, infectious diseases, calcium measurements analysis and pharmacokinetic models, to mention a few (see references above). We addressed both theoretical and practical aspects.
\par
We characterized a necessary and sufficient condition for identifiability of parameters. Specifically, we showed that uniqueness of parameters is not equivalent to uniqueness of ODEs solutions; this fact seems not to have been noticed in previous statistical literature. Exploiting the linearity feature of the model, we developed an integral based two-step estimation approach. The method is based on first estimating the function that is modeled as a solution of the system and then estimating the parameters. It results in an estimator that needs no repeated numerical integration of the system. Moreover, it is consistent at a $\sqrt n$-rate, provided the estimator of the function that solves the system, is sufficiently accurate.

We have studied two specific, sufficiently accurate estimators of the solution of the system, namely a local polynomial estimator (smooth estimator) and an estimator based on averages (step function estimator). We call our estimators "modified integral methods". Although the size of the system in terms of the dimensions $d$ and $p$ does not matter in the theoretical results, in practice it makes a difference, since computing time will grow with these dimensions. However, this growth will be modest since our "modified integral methods" do not employ search algorithms. We studied both estimation approaches via numerical simulations. We compared the smooth estimator to the derivative based two-step approach and to the generalized profiling method. The finite sample performance of the integral estimator is substantially better than that of the derivative based method. As expected, the variability of the generalized profiling approach is smaller; however, it requires a complex optimization step that can affect the estimation results if started too far from the 'true' vector of parameters. Therefore it makes sense to use the integral estimator in order to generate a preliminary estimator to be used as an initial guess for the optimization step of other, more complicated, but accurate,  estimation approaches. The step function estimator was tested under several scenarios of experimental studies; the numerical results support the theory and suggest that the estimation accuracy is robust with respect to the distribution of the errors. Furthermore, we see that practically, the number of repeated measures may be relatively small without the accuracy being corrupted. All simulations were executed in Matlab. The code for executing these simulations and for implementing the method for user data is added as supplementary material to this paper.
\section*{Appendix A: Proofs}\label{sec:proofs}
%
%
\subsection*{A.1. Proof of Proposition~\ref{pr:identifiability}}
{\rm (i)}
Let $C_W$ be nonsingular. If $A_W - B_W C_W^{-1}B_W^T$ would be
singular, there would exist a $d$-vector $\zeta \neq 0$ with
$A_W \zeta=B_W C_W^{-1}B_W^T \zeta.$ Consequently, in view of the nonsingularity of $A_W$ we would have $B_W^T \zeta\neq 0$
and $C_WC_W^{-1}B_W^T \zeta=B_W^T A_W^{-1} B_W C_W^{-1}B_W^T\zeta.$ With
$\eta=C_W^{-1}B_W^T \zeta$ and because $C_W$ is nonsingular, this implies
$\eta \neq 0$ and $C_W\eta=B_W^T A_W^{-1} B_W \eta.$ Consequently, we obtain
\begin{eqnarray}\label{Schurrelation}
\lefteqn{0=\eta^T\left(C_W-B_W^T A_W^{-1}B_W \right)\eta} \nonumber \\
&&=\int_0^1 \eta^T \left(G(t)- A_W^{-1} B_W \right)^T \,\rd W(t) \left(G(t)- A_W^{-1} B_W \right)\eta \nonumber \\
&& =\ \parallel (G - A_W^{-1} B_W) \eta \parallel_W^2
\end{eqnarray}
in the notation of (\ref{innerproduct}) and in view of (\ref{normingW}) and (\ref{norm}). Hence the $i$-th component of $\left(G(t)- A_W^{-1} B_W\right)\eta$ vanishes for $W_{ii}$-almost all $t \in [0,1].$ So, the $i$-th component of $G(t)\eta=\int_0^t g(x(s))\,\rd s\,\eta$
is constant and thus equals 0 for $W_{ii}$-almost all $t\in[0,1],$ since 0 belongs to the support of $W_{ii}.$ Because this holds for all $i=1, \dots, d,$ it follows that
\begin{equation*}
C_W \eta=\int_0^1 G^T(t)\,\rd W(t) G(t)\eta=0,
\end{equation*}
which contradicts the nonsingularity of $C_W.$

Equalities (\ref{xi}) and (\ref{theta}) may be verified now by
substituting the right hand side of (\ref{eq:int_model}) for
$x(t).$

{\rm (ii)}
If $C_W$ would be singular, there would exist a
$p$-vector $\eta$ with
\begin{equation*}
\eta^T C_W \eta=\int_0^1 \eta^T G^T(t)\,\rd W(t) G(t)\eta = \parallel G(\cdot) \eta \parallel_W =0,
\end{equation*}
which would imply that the $i$-th component of $G(t)\eta$ vanishes for $W_{ii}$-almost all $t \in[0,1].$ Consequently,
(\ref{eq:int_model}) yields for all $\alpha \in \rR$ and all $i=1,\dots, d$ that the $i$-th component of the equation
\begin{equation*}
x(t)=\xi + \int_0^t g(x(s))\, \rd s (\theta + \alpha \eta)
\end{equation*}
holds for $W_{ii}$-almost all $t \in [0,1],$  and hence knowledge for all $i=1,\dots,d$ of $x_i(t)$ for $W_{ii}$-almost all $t \in [0,1]$ would not determine
$\theta \in \Theta,$ since $\Theta$ is open.
\epr

\begin{remark}
Interestingly, the start of the proof of
Proposition~\ref{pr:identifiability} may also be formulated via
the concept of Schur complement. Let
\[
M=\begin{pmatrix} I_d& B_W\\B_W^T&C_W
\end{pmatrix},
\]
where the entries of the matrix $M$ are defined in (\ref{BC}). If
$C_W$ is nonsingular, then the {\it Schur complement} of $C_W$ with
respect to $M$ is $I_d - B_W C_W^{-1}B_W^T$
(\cite{haynsworth1968schur}). Moreover, note that
\[
\begin{pmatrix}
I_d& B_W\\B_W^T&C_W
\end{pmatrix}
=\begin{pmatrix}
I_d& B_WC_W^{-1}\\0&I_d
\end{pmatrix}
\begin{pmatrix}
I_d - B_W C_W^{-1}B_W^T& 0\\0&C_W
\end{pmatrix}
\begin{pmatrix}
I_d& 0\\C_W^{-1}B_W^T&I_d
\end{pmatrix}.
\]
Taking determinants of both sides it is immediately clear that
${\rm det}(M)={\rm det}(C_W)\times{\rm det}(I_d - B_W C_W^{-1}B_W^T)$.
Consequently, if $C_W$ is nonsingular and $I_d - B_W C_W^{-1}B_W^T$ is
singular then $M$ is singular. This implies that we can find a
vector $(x,y)\ne 0$ such that $I_dx+B_Wy=0$ and $B_W^Tx+C_Wy=0$ (note
that $y\ne 0$ otherwise $x=0$). Solving the first equation for $x$
and plugging into the second equation we obtain
(\ref{Schurrelation}) with $\eta=y.$
\end{remark}
\subsection*{A.2. Proof of Theorem~\ref{th:const}}
Denote the supnorm of $x(\cdot)$ by $M=\sup_{t\in [0,1]}
\parallel x(t)\parallel=\parallel x \parallel_\infty <\infty.$
Since the map $g(\cdot)$ is continuous on $\rR^d,$ it is
continuous on the compact ball $B_{M+1}=\{x\in \rR^d\,|\ \parallel
x\parallel \leq M+1\}.$ Consequently, (each component of)
$g(\cdot)$ is bounded and uniformly continuous on $B_{M+1}.$
\par
Fix $\varepsilon>0.$ There exists a $\delta>0$ such that for all $x,y
\in B_{M+1}$ with $\parallel x-y\parallel <\delta$ the inequality
$\parallel g(x)-g(y)\parallel < \varepsilon$ holds, with the norm
$\parallel \cdot \parallel$ of a matrix equal to the square root
of the sum of squares of the components of the matrix.
Consequently, $\parallel \hat{x}_n-x\parallel_\infty <\delta$
implies $\int_0^1\parallel g(\hat{x}_n(t))-g(x(t))\parallel \,\rd t <
\varepsilon,$ and hence we have
\begin{equation*}\label{uniform}
P\left(\int_0^1\parallel g(\hat{x}_n(t))-g(x(t))\parallel \,\rd t
\geq \varepsilon \right) \leq P\left(\parallel
\hat{x}_n-x\parallel_\infty \geq \delta \right).
\end{equation*}
Together with the consistency (\ref{supnorm}) of
$\hat{x}_n(\cdot)$ this implies
\begin{equation}\label{Ghatconvergence}
\sup_{t\in[0,1]}\parallel \hat{G}_n(t)-G(t)\parallel\ \stackrel
P{\to} 0.
\end{equation}
Since $g(\cdot)$ is bounded on $B_{M+1}$ and $x(\cdot)$ is bounded
on $[0,1],$ so is $G(\cdot).$ Consequently,
(\ref{Ghatconvergence}) yields boundedness of $\hat{G}_n(\cdot)$
on $[0,1]$ in probability. Using the boundedness and continuity of $G(\cdot)$, the boundedness in probability of $\hat{G}_n(\cdot),$ and (\ref{Ghatconvergence}), and applying the weak convergence of $(W_n)$ and dominated convergence we obtain
\begin{equation*}\label{convBC}
\hat{B}_n \stackrel P{\to} B_W,\quad \hat{C}_n \stackrel P{\to} C_W.
\end{equation*}

Since the consistency (\ref{supnorm}) of $\hat{x}_n(\cdot)$ also
implies
$P\left( \parallel \hat{x}_n\parallel_\infty  > M+1 \right) \to
0,$ again by the weak convergence of $(W_n)$ and dominated convergence we obtain the consistency of
(\ref{xihat}) and (\ref{thetahat}).

Let $c>0$ and $C<\infty$ be such that for all $x,y \in {\mathbb R}^d$ the inequalities $c \leq d_n(x,y)/\parallel x-y \parallel \leq C$ hold.
By the triangle inequality for $d_n(\cdot,\cdot)$ and (\ref{nuhat}) we have
\begin{equation}\label{triangle}
d_n(h(\hat{\nu}_n), h(\nu)) \leq d_n(h(\hat{\nu}_n), \hat{\theta}_n) + d_n(\hat{\theta}_n, h(\nu)) \leq 2 d_n(\hat{\theta}_n, \theta) + \frac 1n,
\end{equation}
and hence
\begin{equation}\label{dnEuclid}
\parallel h(\hat{\nu}_n) - h(\nu) \parallel \leq \frac{2C}c \parallel \hat{\theta}_n - \theta \parallel + \frac 1{cn} \stackrel P{\to} 0,
\end{equation}
as $n \to \infty.$ Since $h^{-1}(\cdot)$ is continuous, this implies consistency of $\hat{\nu}_n.$
\epr

\subsection*{A.3. Proof of Theorem~\ref{th:rootnconsistency}}
First, we collect some properties of the semidefinite inner product (\ref{innerproduct}) that we need.
\begin{lemma}\label{CauchySchwarz}
Let $W$ be a symmetric $d \times d$-matrix of finite signed measures on $([0,1], \cal B)$ such that (\ref{innerproduct}) defines a (nonnegative) semidefinite inner product. Then the diagonal elements of $W$ are nonnegative measures on $([0,1], \cal B),$ the Cauchy-Schwarz inequality $|<x,y>_W| \leq \ \parallel x \parallel_W \ \parallel y \parallel_W$ holds, and, in particular, for all $i=1,\dots,d,\ j=1,\dots,d,$ and for all $x_i(\cdot)$ and $x_j(\cdot)$, such that $\int x_i^2 \rd W_{ii}, \int x_j^2 \rd W_{jj},$ and $\int x_i x_j \rd W_{ij}$ are well-defined and finite,
\begin{equation}\label{CS}
\left[\int x_i x_j \rd W_{ij}\right]^2 \leq
\int x_i^2 \rd W_{ii} \int x_j^2 \rd W_{jj}
\end{equation}
holds.
\end{lemma}
\pr
If there would exist a Borel set $B \subset [0,1]$ with $W_{ii}(B)<0,$ then $W_{ii}(B)=<x,x>_W <0$ would hold for $x^T(t) = (0,\dots, 0, {\bf 1}_B(t), 0, \dots, 0)$ with the indicator as the $i$-th component, thus contradicting the nonnegative semidefiniteness of the inner product.

The Cauchy-Schwarz inequality with $x^T(t)=(0,\dots, 0, x_i(t), 0, \dots, 0)$ and $y^T(t)=(0,\dots, 0, x_j(t), 0, \dots, 0)$ reads as (\ref{CS}).
\epr
We continue with another lemma that will be used in the sequel.
\begin{lemma}\label{lem:main_term}
Under the conditions of Theorem \ref{th:rootnconsistency}
\begin{equation}\label{mainterm}
<\hat{G}_n - G, \hat{G}_n - G>_{W_n} = O_p\left(c_n^2 + d_n^2 + v_n^2 \right)
\end{equation}
holds.
\end{lemma}
\pr
The left hand side of (\ref{mainterm}) is a $p\times p$-matrix. Its entry in the $i$-th row and $j$-th column equals
\begin{eqnarray}\label{elementij}
\lefteqn{\sum_{h=1}^d \sum_{k=1}^d \int_0^1 \int_0^t \left[g_{hi}(\hat{x}_n(s)) - g_{hi}(x(s))\right]\rd s}\\
 && \qquad\qquad \int_0^t \left[g_{kj}(\hat{x}_n(u)) - g_{kj}(x(u))\right]\rd u\ \rd W_{n,hk}(t) \nonumber
\end{eqnarray}
In view of the Cauchy-Schwarz inequality (\ref{CS}) this shows that it suffices to prove
\begin{equation}\label{elementii}
\int_0^1 \left\{\int_0^t \left[g_{hi}(\hat{x}_n(s)) - g_{hi}(x(s))\right]\rd s \right\}^2 \rd W_{n,hh}(t)= O_p\left(c_n^2 + d_n^2 + v_n^2 \right).
\end{equation}
Denote by $\partial g_{hi}(x)/\partial x$ the $d$-dimensional
row vector of first derivatives of the entry $g_{hi}$ of the
matrix $g,$ and by $\partial^2 g_{hi}(x)/\partial x^2$ the
$d\times d$-matrix of second derivatives. The following Taylor
expansion holds
\begin{eqnarray}\label{Taylor}
\lefteqn{g_{hi}(\hat{x}_n)=g_{hi}(x)+\frac \partial{\partial x}
g_{hi}(x)\left( \hat{x}_n -x\right)}\\ \nonumber && +\int_0^1
\int_0^\lambda \left( \hat{x}_n -x\right)^T
\frac{\partial^2}{\partial x^2}g_{hi}\left(x + \mu\left(\hat{x}_n
-x\right)\right) \left( \hat{x}_n -x\right)\rd \mu\, \rd \lambda.
\end{eqnarray}
In view of $\parallel x\parallel _\infty < \infty$ and the
continuity of the partial derivatives of $g_{hi}(\cdot)$ the
function $s \mapsto
\partial g_{hi}(x(s))/\partial x$ is bounded on [0,1]. By
(\ref{variance}) and (\ref{supnormc}) this implies
\begin{equation}\label{firstderivative}
\int_0^1 \rE\left(\left\{\int_0^t \frac \partial{\partial x}
g_{hi}(x(s))\left( \hat{x}_n(s) -x(s)\right) \rd s \right\}^2 \right) \rd W_{n,hh}(t) =
O\left( v_n^2 + c_n^2 \right).
\end{equation}
Similarly, in view of $\parallel x\parallel _\infty < \infty,$ of
$\parallel \hat{x}_n\parallel _\infty =O_p(1),$ and of the
continuity of the partial second derivatives of $g_{hi}(\cdot)$
the function $s \mapsto
\partial^2 g_{hi}\left(x(s)+\mu\left(\hat{x}_n(s)
-x(s)\right)\right)/\partial x^2$ is bounded in probability. By
(\ref{dn}) and (\ref{supnormc}) this implies
\begin{eqnarray}\label{secondderivative}
\sup_{0\leq t \leq 1}\left|\int_0^t \int_0^1 \int_0^\lambda (
\hat{x}_n(s) -x(s))^T \frac{\partial^2}{\partial
x^2}g_{hi} (x(s) + \mu (\hat{x}_n(s) -x(s))) \right. && \nonumber \\
\left. (\hat{x}_n(s)-x(s))\rd\mu\, \rd \lambda\, \rd s
\right|=O_p\left(d_n + c_n^2 \right). &&
\end{eqnarray}
Combining (\ref{Taylor}), (\ref{firstderivative}), and
(\ref{secondderivative}) we arrive at (\ref{elementii}) and hence (\ref{mainterm}).
Note that $\parallel \hat{x}_n\parallel _\infty =O_p(1)$ is used
in the argument leading up to (\ref{secondderivative}) in order to
obtain boundedness of  $s \mapsto
\partial^2 g_{hi}\left(x(s)+\mu\left(\hat{x}_n(s)
-x(s)\right)\right)/\partial x^2$ in probability. If all second
partial derivatives of all $g_{hi}(\cdot)$ are bounded,
$\parallel \hat{x}_n\parallel _\infty =O_p(1)$ is not needed for
this. Moreover, continuity of the second derivatives is not needed
for this either. The proof of the Lemma is complete. \epr
\paragraph{Proof of Theorem~\ref{th:rootnconsistency}}
We write
\begin{equation}\label{differenceA}
\hat{B}_n - B_W = <I_d, \hat{G}_n - G >_{W_n} + \left\{ <I_d, G>_{W_n} - <I_d,G>_W \right\}.
\end{equation}
Since $g(\cdot)$ is continuous and $x(\cdot)$ is bounded on $[0,1],$ we may conclude that $G(\cdot)$ is differentiable with bounded derivatives. This implies that the second term at the right hand side of (\ref{differenceA}) is of order $O(w_n)$ in view of (\ref{orderweakconvergence}). Each matrix entry of the first term is a sum of $d$ terms of the type
\begin{equation}\label{firstAterm}
\int_0^1 \left(\hat{G}_{n,hj} - G_{n,hj}\right) \rd W_{n,ih}.
\end{equation}
By the Cauchy-Schwarz inequality of Lemma \ref{CauchySchwarz} and by Lemma \ref{lem:main_term}, in particular formula (\ref{elementii}), we see that each such a term is of order $O_p\left(c_n + d_n +v_n\right).$ We have shown
\begin{equation}\label{orderdifferenceA}
\hat{B}_n -B_W = O_p\left(c_n + d_n + v_n + w_n \right).
\end{equation}
Similarly we study
\begin{eqnarray}\label{differenceB}
\lefteqn{\hat{C}_n - C_W = <\hat{G}_n - G, \hat{G}_n - G>_{W_n}} \\
&& \quad \qquad +2<G,\hat{G}_n - G>_{W_n} + \left\{<G,G>_{W_n} - <G,G>_W \right\}.\nonumber
\end{eqnarray}
The second and third term at the right hand side are handled by the same arguments as the first and second term at the right hand side of (\ref{differenceA}), respectively. Consequently, by Lemma \ref{lem:main_term} we arrive at
\begin{equation}\label{orderdifferenceB}
\hat{C}_n -C_W = O_p\left(c_n + d_n + v_n +w_n \right).
\end{equation}
We also study
\begin{eqnarray}\label{Gx}
\lefteqn{<\hat{G}_n, \hat{x}_n>_{W_n} - <G,x>_W = <\hat{G}_n - G, \hat{x}_n>_{W_n}} \\
&& \quad \qquad + <G,\hat{x}_n - x>_{W_n} + \left\{<G,x>_{W_n} - <G,x>_W \right\}.\nonumber
\end{eqnarray}
By (\ref{supnormc}), Lemma \ref{CauchySchwarz}, and Lemma \ref{lem:main_term} the first term at the right hand side of (\ref{Gx}) is of the order $O_p(c_n + d_n + v_n).$ The $i$-th component of the $p$-vector that is the second term, is a sum of $d^2$ terms of the type
\begin{equation}\label{replacexhat}
\int_0^1 G_{hi}(t)\left[\hat{x}_{nj}(t) -x_j(t)\right] \rd W_{n,hj}(t).
\end{equation}
Since $G(\cdot)$ is continuous and hence bounded on $[0,1],$ we obtain by (\ref{supnormc}) and (\ref{variance2}), that (\ref{replacexhat}) and hence the second term at the right hand side of (\ref{Gx}) is of order $O_p(c_n + v_n).$  The third term at the right hand side of (\ref{Gx}) is of order $O(w_n)$ in view of (\ref{orderweakconvergence}), where we note that both $G(\cdot)$ and $x(\cdot)$ are differentiable with bounded derivatives. We have obtained
\begin{equation}\label{orderdifferenceGx}
<\hat{G}_n, \hat{x}_n>_{W_n} - <G,x>_W = O_p\left(c_n + d_n + v_n +w_n \right).
\end{equation}
In a similar way we obtain
\begin{equation}\label{orderdifferenceIx}
<I_d, \hat{x}_n>_{W_n} - <I_d,x>_W = O_p\left(c_n + v_n +w_n \right).
\end{equation}

Writing $\hat \theta_n -\theta$ and $\hat \xi_n - \xi$ as
telescoping sums in which sequentially random elements are
replaced by the corresponding deterministic ones, and applying
(\ref{orderdifferenceA}), (\ref{orderdifferenceB}), (\ref{orderdifferenceGx}), and (\ref{orderdifferenceIx}) repeatedly, we obtain a proof of the consistency to the order
$O_p\left(c_n + d_n + v_n + w_n \right)$ of $\hat{\theta}_n$ and $\hat{\xi}_n.$ Subsequently the consistency of
$\hat{\nu}_n$ to the same order is obtained via (\ref{dnEuclid}) and the Lipschitz continuity of $h^{-1}(\cdot).$
\epr
\subsection*{A.4. Proof of Theorem~\ref{cor:local}}

The following lemma assures us that the local polynomial
estimator $\hat{x}_n$ satisfies the conditions as required in
Theorem~\ref{th:rootnconsistency}.
\begin{lemma}\label{lem:nonpar}
Let the model be defined by \rf{ode_model}--\rf{int_model}.
Suppose that for any $j=1,...,d$ the solution $x_j(t;\theta,\xi)$
is a $C^\alpha$-function of $t$ on the interval $[0,1]$ for some
positive real $\alpha \geq 1$.
\par
Let $W_n$ be as in Theorem \ref{cor:local}, let the observations be given by (\ref{eq:observations}) where we
have $t_i=i/n$, $i=1,...,n,$ and let the estimator for $x(\cdot)$ be
defined in \rf{est_x}. Assume that the errors
$\varepsilon_j(t_i),\ i=1,...,n,\ j=1,...,d,$ are i.i.d. and have
mean $\rE\varepsilon_j(t_i)=0$ and variance
$\rE\varepsilon_j(t_i)^2=\sigma^2_\varepsilon<\infty$. Under
Condition~$K$
\begin{equation}\label{supnormxhat}
\parallel \hat{x}_n -x \parallel _\infty =O_p\left(\frac1{n^{1/3}b}\right),
\end{equation}
\begin{equation}\label{cn}
\parallel \rE\hat{x}_n - x \parallel_\infty =O(b^\alpha),
\end{equation}
and
\begin{equation}\label{dn2}
\rE \left(\parallel\hat{x}_n(t) -x(t)\parallel_{W_n}^2
\right) = O\left(b^{2\alpha}+\frac 1{nb}\right)
\end{equation}
hold, and for any $h, j=1,...,d$ and any bounded measurable function
$f(\cdot)$
\begin{equation}\label{varf}
\int_0^1{\rm var}\Big(\int_0^t f(s)\hat{x}_{n,j}(s)\rd s\Big)\rd W_{n,hh}(t)
=O\left(\frac 1n \right)
\end{equation}
and
\begin{equation}\label{varf2}
{\rm var}\Big(\int_0^1 f(t)\hat{x}_{n,j}(t)\rd W_{n,hj}(t) \Big) =
O\left(\frac 1n \right)
\end{equation}
hold.
\end{lemma}
\pr Our proof is based on \cite[Chapter 1]{tsybakov2009introduction}. In
particular, the proofs of (\ref{cn}) and (\ref{dn2}) follow from
his Proposition 1.13. Note that the bounds given in this Proposition 1.13 are uniform over $[0,1],$ and that (\ref{dn2}) needs an application of Fubini's theorem and the boundedness of $\int_0^1 \rd W_{n,hh}(t)$ as guaranteed by (\ref{Wrootn1}) and the finiteness of the entries of $W.$

Lemma 1.5 and (1.70) of \cite{tsybakov2009introduction} show that
Condition~{\rm K(iii)} implies that there exists a positive
integer $n_0$ and a positive constant $\lambda_0$ such that for
all $n\geq n_0$, $t\in[0,1]$ and $v\in\rR^{\ell+1}$ the inequality
$\parallel B_n(t)^{-1}v\parallel\leq\parallel
v\parallel\lambda_0^{-1}$ holds, where $\parallel\cdot\parallel$
stands for the Euclidean norm in $\rR^{\ell+1}$. This together
with $\parallel U(0)\parallel=1$ leads for $n\geq n_0$ to
\begin{eqnarray*}
|V_{n,i}(t)|&=&\Big|
\frac{1}{nb}U^T(0)B_n(t)^{-1}U\Big(\frac{t_i-t}{b}\Big)K\Big(\frac{t_i-t}{b}\Big)\Big|
\\&\leq&\frac{1}{nb\lambda_0}\parallel
U\Big(\frac{t_i-t}{b}\Big)K\Big(\frac{t_i-t}{b}\Big)
\parallel
\\&=&
\frac{1}{nb\lambda_0}\parallel
U\Big(\frac{t_i-t}{b}\Big)K\Big(\frac{t_i-t}{b}\Big)\ind_{[|(t_i-t)/b|\leq 1]}
\parallel,
\end{eqnarray*}
since the kernel $K$ is supported on $[-1,1]$. Furthermore,
\[\parallel U\Big(\frac{t_i-t}{b}\Big)\ind_{[|(t_i-t)/b|\leq 1]}
\parallel^2 \leq\sum_{k=0}^\ell\frac{1}{(k!)^2}\leq e\]
holds and hence for sufficiently large $n$
\begin{eqnarray*}
|V_{n,i}(t)|&\leq&\frac{\sqrt{e}}{nb\lambda_0}\Big|K\Big(\frac{t_i-t}{b}\Big)\Big|.
\end{eqnarray*}
Using this bound for $|V_{n,i}(t)|$ we obtain
\begin{eqnarray*}
\lefteqn{\int_0^1{\rm var}\Big(\int_0^t f(s)\hat{x}_{n,j}(s)\rd
s\Big)\rd W_{n,hh}(t)} \\
&& = \int_0^1{\rm var}\Big(\int_0^t
f(s)\sum_{i=1}^nY_j(t_i)V_{n,i}(s) \rd s\Big)\rd W_{n,hh}(t) \\
&& \leq \int_0^1\frac{e
\sigma_\varepsilon^2}{n^2\lambda_0^2}\sum_{i=1}^n \Big(\int_0^t
|f(s)|\frac{1}{b}\Big|K\Big(\frac{t_i-s}{b}\Big)\Big|\rd
s\Big)^2\rd W_{n,hh}(t) \\
&& \leq \frac{e \sigma_\varepsilon^2}{n\lambda_0^2}
\parallel f\parallel^2_\infty\Big(\int_0^1 |K(u)|\rd u\Big)^2 \int_0^1 \rd W_{n,hh}(t)
\end{eqnarray*}
and hence (\ref{varf}) by boundedness of $f,$ by finiteness of
$\sigma_\varepsilon^2$ and $\int_0^1 |K(u)|\rd u,$ which is
implied by Condition~${\rm K(iii)},$ and by boundedness of $\int_0^1 \rd W_{n,hh}(t)$ as guaranteed by (\ref{Wrootn1}) and the finiteness of the entries of $W.$ Similarly and by the Cauchy-Schwarz inequality (\ref{CS}) we obtain
\begin{eqnarray*}
\lefteqn{{\rm var}\Big(\int_0^1 f(t)\hat{x}_{n,j}(t)\rd W_{n,hj}(t) \Big) =
{\rm var}\Big(\int_0^1 f(t)\sum_{i=1}^nY_j(t_i)V_{n,i}(t) \rd W_{n,hj}(t) \Big)} \\
&& = \sum_{i=1}^n E\left( \int_0^1 f(t)\varepsilon_j(t_i)V_{n,i}(t) \rd W_{n,hj}(t) \right)^2 \\
&& \leq \sum_{i=1}^n E\left( \int_0^1 f^2(t)\varepsilon_j^2(t_i) \left|V_{n,i}(t)\right| \rd W_{n,jj}(t) \int_0^1 \left|V_{n,i}(t)\right| \rd W_{n,hh}(t) \right)\\
&&\leq \frac{e
\sigma_\varepsilon^2}{n^2\lambda_0^2} \parallel f\parallel^2_\infty \sum_{i=1}^n \int_0^1
 \frac{1}{b}\left|K\left(\frac{t_i-t}{b}\right) \right| \rd W_{n,jj}(t) \\
 && \qquad \qquad \qquad \qquad \qquad  \int_0^1 \frac{1}{b}\left|K\left(\frac{t_i-t}{b}\right) \right| \rd W_{n,hh}(t) \\
&& = O(n^{-1}),
\end{eqnarray*}
where the last equality holds in view of (\ref{Wrootn2}) since $t \mapsto K((t_i -t)/b)$ is bounded and vanishes outside an interval of length at most $2b.$ We have proved (\ref{varf2}).

To prove (\ref{supnormxhat}), we note that Proposition 1.12,
(1.82), and (1.83) of \cite{tsybakov2009introduction} yield
\begin{eqnarray}\label{a}
\lefteqn{\parallel \hat{x}_n -x \parallel _\infty =
\sup_{t\in[0,1]}\parallel \sum_{i=1}^n V_{n,i}(t)
\varepsilon(t_i)\parallel \nonumber}\\
&&=O_p\left(\frac 1{Mb^2} + \max_{1\leq j \leq M} \parallel
\sum_{i=1}^n V_{n,i}\left(\frac jM \right)
\varepsilon(t_i)\parallel \right).
\end{eqnarray}
Proposition 1.12 of ibid. also implies
\begin{eqnarray}\label{b}
\lefteqn{\rE \left( \max_{1\leq j \leq M} \parallel \sum_{i=1}^n
V_{n,i}\left(\frac jM \right) \varepsilon(t_i)\parallel \right)^2
\leq \sum_{j=1}^M \sum_{i=1}^n V_{n,i}^2\left(\frac jM \right) d
\sigma_\varepsilon^2 \nonumber}\\
&& \leq \frac{d\sigma_\varepsilon^2\sqrt{e}K_{\rm
max}}{nb\lambda_0} \sum_{j=1}^M \sum_{i=1}^n| V_{n,i}\left(\frac jM
\right)|=O\left(\frac M{nb}\right).
\end{eqnarray}
Choosing $M=O(n^{1/3}b^{-1})$ we see that (\ref{a}) and (\ref{b})
imply (\ref{supnormxhat}). \epr

To prove Theorem 3 we first note that (\ref{orderweakconvergence}) is satisfied with $w_n = 1/\sqrt n$ in view of (\ref{Wrootn1}).
Applying Theorem 2 we see that Lemma \ref{lem:nonpar} with $b=b_n=n^{-\beta}$
implies Theorem 3, if the following choices are being made:
\begin{enumerate}
\item[Case 1.]
$1/(2\alpha) \leq \beta \leq 1/3,$
\item[Case 2.]
$1/(2\alpha) \leq \beta \leq 1/2.$
\end{enumerate}
The optimal convergence rate for (\ref{dn2}) is
$n^{-2\alpha/(2\alpha+1)},$ which is obtained by $b=b_n=
n^{1/(2\alpha+1)}.$ Compared to this, undersmoothing is needed to
control the bias in (\ref{cn}).
\par
\subsection*{A.5. Proof of Theorem~\ref{th:step}}

To prove this theorem we apply Theorem~\ref{th:rootnconsistency} again.
As in the preceding proof we first note that (\ref{orderweakconvergence})
is satisfied with $w_n = 1/\sqrt n$ in view of (\ref{Wrootn1}).
Since $g(\cdot)$ is continuous and $x(\cdot)$ is bounded, $g(x(\cdot))$ is. Consequently, we have
\begin{eqnarray}\label{biasI}
\lefteqn{ \parallel \rE \hat{x}_n -x
\parallel_\infty^2=\max_{1\leq i \leq I}\sup_{(i-1)/I \leq t \leq
i/I} \parallel x(i/I)-x(t)\parallel^2 \nonumber}\\
&& = \max_{1\leq i \leq I}\sup_{(i-1)/I \leq t \leq i/I} \parallel
\int_t^{i/I} g(x(s))\rd s \, \theta\parallel^2 \\
&& \leq \max_{1\leq i \leq I} \frac 1I \int_{(i-1)/I}^{i/I}
\parallel g(x(s))\parallel^2\rd s \parallel \theta \parallel^2
=O\left( \frac 1{I^2}\right)= O\left(\frac 1n \right).\nonumber
\end{eqnarray}
Furthermore, (\ref{samplesizes}) and
\begin{eqnarray*}\label{Pbounded}
\lefteqn{P\left(\parallel \hat{x}_n -\rE \hat{x}_n
\parallel_\infty \geq M \right) =P\left( \max_{1\leq i \leq I}
\parallel \frac 1{J_i} \sum_{j=1}^{J_i} \varepsilon^{(j)}(i/I)
\parallel \geq M \right) \nonumber }\\
&& \qquad \leq 1- \prod_{i=1}^I \left(
1-\frac{d\sigma_\varepsilon^2}{M^2 J_i} \right) \leq
\frac{d\sigma_\varepsilon^2}{M^2} \sum_{i=1}^I \frac 1{J_i}
\end{eqnarray*}
show that $\hat{x}_n(\cdot) -\rE \hat{x}_n(\cdot)$ is bounded in
probability. Together with (\ref{biasI}) and the boundedness of
$x(\cdot)$ this proves that $\hat{x}_n(\cdot)$ is bounded in
probability, i.e.
\begin{equation}\label{xIbdd}
\parallel \hat{x}_n \parallel_\infty = O_p(1).
\end{equation}

By (\ref{Wrootn2}) and (\ref{samplesizes}) we also have
\begin{eqnarray}\label{EI}
\lefteqn{\rE \left( < \hat{x}_n - \rE\hat{x}_n, \hat{x}_n - \rE\hat{x}_n>_{W_n} \right) \nonumber}\\
&& = \sum_{i=1}^I \sum_{h=1}^d \sum_{k=1}^d \int_{(i-1)/I}^{i/I}
\rE\left( \frac 1{J_i} \sum_{j=1}^{J_i} \varepsilon_h^{(j)}(\frac iI)
 \frac 1{J_i} \sum_{\ell=1}^{J_i} \varepsilon_k^{(\ell)}(\frac iI)\right) \rd W_{n,hk}(t) \nonumber\\
 && = \sum_{i=1}^I \frac{\sigma_\varepsilon^2}{J_i} \sum_{h=1}^d \int_{(i-1)/I}^{i/I} \rd W_{n,hh}(t) \nonumber \\
&& \leq \frac{C d\sigma_\varepsilon^2}I \sum_{i=1}^I
\frac 1{J_i} = O\left(\frac 1I\right)=O\left(\frac 1{\sqrt
n}\right).
\end{eqnarray}

For any bounded measurable function $f(\cdot)$ and the $j$th
component $\hat{x}_{n,j}(\cdot)$ of $\hat{x}_n(\cdot)$ we obtain
by (\ref{samplesizes})
\begin{eqnarray}\label{varI}
\lefteqn{\int_0^1 {\rm var}\left( \int_0^t f(s)\hat{x}_{n,j}(s)
\rd s \right) \rd W_{n,hh}(t) \nonumber} \\
&& = \sum_{i=1}^I \int_{(i-1)/I}^{i/I}{\rm
var}\left(\sum_{\ell=1}^i \int_{(\ell-1)/I}^{(\ell/I)\wedge t}
 f(s) \rd s \frac 1{J_\ell} \sum_{m=1}^{J_\ell}
\varepsilon_j^{(m)}(\ell/I) \right) \rd W_{n,hh}(t)  \\
&& = \sum_{i=1}^I \int_{(i-1)/I}^{i/I} \sum_{\ell=1}^i \left(
\int_{(\ell-1)/I}^{(\ell/I)\wedge t}
 f(s) \rd s \right)^2 \frac {\sigma_\varepsilon^2}{J_\ell} \rd W_{n,hh}(t) \nonumber \\
&& \leq \frac{C \sigma_\varepsilon^2 \parallel f
 \parallel_\infty^2}{I^2} \sum_{i=1}^I \frac 1{J_i}
=O\left( \frac 1{I^2}\right)= O\left(\frac 1n \right) \nonumber
\end{eqnarray}
as well. Similarly we get
\begin{eqnarray}\label{varI2}
\lefteqn{ {\rm var}\left( \int_0^1 f(t)\hat{x}_{n,j}(t) \rd W_{n,hj}(t)
\right) \nonumber}\\
&& = {\rm var}\left(\sum_{i=1}^I \int_{(i-1)/I}^{i/I}
 f(t) \rd W_{n,hj}(t) \frac 1{J_i} \sum_{m=1}^{J_i}
\varepsilon_j^{(m)}(i/I) \right)  \\
&& = \sum_{i=1}^I  \left( \int_{(i-1)/I}^{i/I}
 f(t) \rd W_{n,hj}(t) \right)^2 \frac {\sigma_\varepsilon^2}{J_i} \nonumber \\
&& \leq \frac{C^2 \sigma_\varepsilon^2 \parallel f
 \parallel_\infty^2}{I^2} \sum_{i=1}^I \frac 1{J_i}
  =O\left( \frac 1{I^2}\right)= O\left(\frac 1n \right). \nonumber
\end{eqnarray}
Applying Theorem~\ref{th:rootnconsistency} with $c_n=O(1/\sqrt n), d_n=
O(1/\sqrt n), v_n=O(1/\sqrt n),$ and $w_n = O(1/\sqrt n),$ we arrive at
Theorem~\ref{th:step}. \epr

\section*{Acknowledgements}
This research was supported by the Dutch Technology Foundation STW, which is part of the Netherlands Organisation for Scientific Research (NWO) and which is partly funded by the Ministry of Economic Affairs.

This research started when the first author was a Postdoc at EURANDOM, Eindhoven University of Technology, and the second one was a Senior Fellow there.

\bibliographystyle{chicago}
\bibliography{bib}

\end{document}